\documentclass[preprint]{imsart}

\RequirePackage[OT1]{fontenc}
\RequirePackage{amsfonts}
\RequirePackage[numbers]{natbib}
\RequirePackage[colorlinks,citecolor=blue,urlcolor=blue]{hyperref}
\RequirePackage[centertags, leqno, sumlimits, nointlimits, namelimits]{mathtools}
\RequirePackage[thmmarks,thref,amsthm,amsmath]{ntheorem}
\RequirePackage{paralist}
\RequirePackage[english, noabbrev]{cleveref}
\usepackage{pdfsync} 
\arxiv{}

\startlocaldefs
\numberwithin{equation}{section}
\theoremstyle{plain}
\newtheorem{theorem}{Theorem}[section]
\newtheorem{lemma}[theorem]{Lemma}
\newtheorem{corollary}[theorem]{Corollary}

\theoremstyle{remark}
\newtheorem{remark}[theorem]{Remark}

\newtheorem{example}[theorem]{Example}

\endlocaldefs

\newcommand{\N}{\mathbb{N}} %Natural numbers
\newcommand{\R}{\mathbb{R}} %Real numbers
\newcommand{\E}{{\mathbb{E}}} %Expectation
\begin{document}

\begin{frontmatter}
\title{On the consistency of adaptive multiple tests}
\runtitle{Consistency of adaptive multiple tests}

\begin{aug}
\author{\fnms{Marc} \snm{Ditzhaus}\ead[label=e1]{marc.ditzhaus@uni-duesseldorf.de}}
\and
\author{\fnms{Arnold} \snm{Janssen}\ead[label=e2]{janssena@uni-duesseldorf.de}}

\runauthor{M. Ditzhaus and A. Janssen}

\affiliation{Heinrich-Heine University D\"usseldorf \thanksmark{m1}}

\address{M. Ditzhaus and A. Janssen\\
Heinrich-Heine University D\"usseldorf\\
Mathematical Institute\\
Universit\"atststra{\ss}e 1\\
40225 D\"usseldorf, Germany\\
\printead{e1}\\
\phantom{E-mail:\ }\printead*{e2}
}

\end{aug}

\begin{abstract}
	Much effort has been done to control the "false discovery rate" (FDR) when $m$ hypotheses are tested simultaneously. The FDR is the expectation of the "false discovery proportion"  $\text{FDP}=V/R$ given by the ratio of the number of false rejections $V$ and all rejections $R$. In this paper, we have a closer look at the FDP for adaptive linear step-up multiple tests. These tests extend the well known Benjamini and Hochberg test by estimating the unknown amount $m_0$ of the true null hypotheses. We give exact finite sample formulas for higher moments of the FDP and, in particular, for its variance. Using these allows us a precise discussion about the consistency of adaptive step-up tests. We present sufficient and necessary conditions for consistency on the estimators $\widehat m_0$ and the underlying probability regime. We apply our results to convex combinations of generalized Storey type estimators with various tuning parameters and (possibly) data-driven weights. The corresponding step-up tests allow a flexible adaptation. Moreover, these tests control the FDR at finite sample size. We compare these tests to the classical  Benjamini and Hochberg test and discuss the advantages of it. 
\end{abstract}

\begin{keyword}[class=MSC]
\kwd[Primary ]{62G10}
\kwd[; secondary ]{62G20}
\end{keyword}

\begin{keyword}
\kwd{false discovery rate (FDR)} 
\kwd{false discovery proportion (FDP)} 
\kwd{multiple testing}
\kwd{adaptive Benjamini Hochberg methods}
\kwd{Storey estimator}
\kwd{$p$-values}
\end{keyword}

\end{frontmatter}

\section{Introduction}\label{sec:intro}
Testing $m\geq 2$ hypotheses simultaneously is a frequent issue in statistical practice, e.g. in genomic research. A widely used criterion for deciding which of these hypotheses should be rejected is the so-called "false discovery rate" (FDR) promoted by \citet{Benjamini_Hochberg_1995}. The FDR is the expectation of the "false discovery proportion" (FDP), the ratio $\text{FDP}=V/R$ of the number of false rejections $V$ and all rejections $R$. Let a level $\alpha\in(0,1)$ be given. Under the so-called basic independence (BI) assumption we have FDR$=(m_0/m)\alpha$ for the classical Benjamini and Hochberg linear step-up test, briefly denoted by BH test. Here, $m_0$ is the unknown amount of true null hypotheses. To achieve a higher power it is of great interest to exhaust the FDR as good as possible. Especially, if $m_0/m=c$ is not close to $1$ there is space for improvement of the BH test. That is why since the beginning of this century the interest of adaptive tests grows. The idea is to estimate $m_0$ by an appropriate estimator $\widehat m_0$ in a first step and to apply the BH test for the (data depended) level $\alpha'=(m/\widehat m_0)\alpha$ in the second step. Heuristically, we obtain for a good estimator $\widehat m_0\approx m_0$ that FDR$\;\approx \alpha$. \citet{Benjamini_Hochberg_2000} suggested an early estimator for $m_0$ leading to an FDR controlling test. Before them \citet{SchwederSpjotvoll1982} already discussed estimators for $m_0$ using plots of the empirical distribution function of the $p$-values. The number of estimators suggested in the literature is huge, here only a short selection: \citet{BenjaminiETAL2006}, Blanchard and Roquain \cite{BlanchardRoquain2008,BlanchardRoquain2009} and \citet{ZeiselETAL2011}. We want to emphasize the papers of the \citet{Storey2002} and \citet{StoreyETAL2004} and, in particular, the Storey estimator based on a tuning parameter $\lambda$. We refer to \citet{StoreyTibshirani2003} for a discussion of the adjustment of the tuning parameter $\lambda$. Generalized Storey estimators with data dependent weights, which were already discussed by \citet{HeesenJanssen2016}, will be our prime example for our general results. A nice property of them is that we have finite sample FDR control, see \cite{HeesenJanssen2016}. Sufficient conditions for finite sample FDR control on general estimators $\widehat m_0$ can be found in \citet{Sarkar2008} and Heesen and Janssen \cite{HeesenJanssen2015,HeesenJanssen2016}. \\
Beside the FDR control there are also other control criteria, for example  the family-wise error rate FWER$=P(V>0)$. Also for the control of FWER adaptive tests, i.e. tests using an (plug-in) estimator for $m_0$, are used and discussed in the literature, see e.g. \citet{FinnerGontscharuk2009} and \citet{SarkarETAL2012}. \\
Stochastic process methods were applied to study the asymptotic behavior of the FDP, among others to calculate asymptotic confidence intervals, and the familywise error rate (FWER) in detail, see \citet{GenoveseWassermann2004}, \citet{MeinshausenBuehlmann2005}, \citet{MeinshausenRice} and Neuvial \cite{Neuvial2008}. Dealing with a huge amount of $p$-values the fluctuation of the FDP becomes, of course, relevant. \citet{FerreiraZwinderman2006} presented formulas for higher moments of FDP for the BH test and \citet{RoquainVillers2011} did so for step-up and step-down tests with general (but data independent) critical values. We generalize these formulas to adaptive step-up tests using general estimators $\widehat m_0$ for $m_0$. In particular, we derive an exact finite sample formula for the variability of FDP. As an application of this we discuss the consistency of FDP and present sufficient and necessary conditions for it. We also discuss the more challenging case of sparsity in the sense that $m_0/m\to 1$ as $m\to \infty$. This situation can be compared to the one of \citet{AbramovichETAL2006}, who derived an estimator of the (sparse) mean of a multivariate normal distribution using FDR procedures.

\textbf{Outline of the results}. In \Cref{sec:preliminaries} we introduce the model as well as the adaptive step-up tests, and in particular the generalized Storey estimators which serve as prime examples. \Cref{sec:moments} provides exact finite sample variance formulas for the FDP under the BI model. Extensions to higher moments can be found in the appendix, see \Cref{sec:higher_moments}. These results apply to the variability and the consistency of FDP, see Section \Cref{sec:consis}. Roughly speaking we have consistency if we have stable estimators $\widehat m_0/m\approx C_0$ and the number of rejections tends to infinity. \Cref{sec:estimators} is devoted to concrete adaptive step-up tests mainly based on the convex combinations of generalized Storey estimators with data dependent weights. We will see that consistency cannot be achieved in general. Under mild assumptions the adaptive tests based on the estimators mentioned above are superior compared to the BH test: The FDR is more exhausted but remains finitely controlled by the level $\alpha$. Furthermore, they are consistent at least when the BH test is consistent. In \Cref{sec:LFC} we discuss least favorable configurations which serve as useful technical tools. For the reader's convenience we add a discussion and summary of the paper in \Cref{sec:discussion}. All proofs are collected in \Cref{sec:proof}.

\section{Preliminaries}\label{sec:preliminaries}
\subsection{The model and general step-up tests}
Let us first describe the model and the procedures. A multiple testing problem consists on $m$ null hypotheses $(H_{i,m},p_{i,m})$ with associated $p$-value $0\leq p_{i,m}\leq 1$ on a common probability space $(\Omega_m, \mathcal{A}_m, P_m)$. We will always use the basic 
independence \textbf{(BI)} assumption given by
\begin{enumerate}[1]
	\renewcommand\labelenumi{\bfseries(\text{BI}\theenumi)} %Dadurch wird alles dick
	\item\label{BI1} The set of hypotheses can be divided in the disjoint union $I_{0,m}\bigcup I_{1,m}=\{1,...,m\}$ of unknown portions of true null $I_{0,m}$ and false null $I_{1,m}$, respectively. Denote by $m_{j}=\#I_{j,m}$ the cardinality of $I_{j,m}$ for $j=0,1$.	
	
	\item\label{BI2} The vectors of $p$-values $\left(p_{i,m}\right)_{i\in I_{0,m}}$ and $\left(p_{i,m}\right)_{i\in I_{1,m}}$ are independent, where each dependence structure is allowed for the $p$-values $\left(p_{i,m}\right)_{i\in I_1}$ for the false hypotheses.
	
	\item\label{BI3} The $p$-values $\left(p_{i,m}\right)_{i\in I_{0,m}}$ of the true null are independent and uniformly distributed on $[0,1]$,  i.e. $P_m(p_{i,m}\leq x)= x$ for all $x\in[0,1]$.

\end{enumerate}
Throughout the paper let $m_0\geq 1$ be nonrandom. As in \citet{HeesenJanssen2016} the results can be extended to more general models with random $m_0$ by conditioning under $m_0$. By using this modification the results easily carry over to familiar mixture models discussed, for instance, by \citet{AbramovichETAL2006} and \citet{GenoveseWassermann2004}. We study adaptive multiple step-up tests with estimated critical values extending the famous \citet{Benjamini_Hochberg_1995} step-up test, briefly denoted by BH test. In the following we recall the definition of this kind of tests. Let
\begin{equation}\label{stepup001}
	0=\alpha_{0:m} < \alpha_{1:m} \leq \alpha_{2:m} \leq \ldots \leq \alpha_{m:m} < 1
\end{equation}
denote possibly data dependent critical values. As an example for the critical values we recall the ones for the BH test, which do not depend on the data:
\begin{align*}
\alpha_{i:m}^{\text{BH}}=\frac{i}{m}\alpha.
\end{align*}
If $p_{1:m}\leq p_{2:m}\leq\ldots\leq p_{m:m}$ denote the ordered $p$-value then the number of rejections is given by
\begin{align*}
	R_m := \max\{ i=0,\ldots,m: p_{i:m}\leq \alpha_{i:m}\}, \text{ where }p_{0:m}:=0,
\end{align*}
and the multiple procedure acts as follows:
\begin{align*}
	\text{reject }H_{i,m}\text{ iff } p_{i,m}\leq \alpha_{R_m:m}.
\end{align*}
Moreover, let 
\begin{align}\label{eqn:def_Vm}
	V_m:= \# \{i\in I_{0,m}\cup\{0\}: p_{i:m}\leq \alpha_{i:m}\}
\end{align}
be the number of falsely rejected null hypothesis. Then the false discovery rate $\text{FDR}_m$ and the false discovery proportion $\text{FDP}_m$ are given by
\begin{align}\label{eqn:def_FDP_FDR}
	&\text{FDP}_m= \frac{V_m}{R_m}\text{ and }\text{FDR}_m=\E\Bigl( \frac{V_m}{R_m}\Bigr)\text{ with }\frac{0}{0}=0.
\end{align}
Good multiple tests like the BH test or the frequently applied adaptive test of \citet{StoreyETAL2004} control the FDR at a pre-specified acceptance error bound $\alpha$ at least under the BI assumption. Besides the control, two further aspects are of importance and discussed below:
\begin{enumerate}[(i)]
	\item To make the test sensitive for signal detection the FDR should exhaust the level $\alpha$ as best as possible.
	
	\item On the other hand the variability of the FDP, see \eqref{eqn:def_FDP_FDR}, is of interest in order to judge the stability of the test. 
\end{enumerate}
For a large class of adaptive tests exact FDR formulas were established in \citet{HeesenJanssen2016}. Here the reader can find new $\alpha$-controlling tests with high FDR. These formulas are now completed by formulas for exact higher FDP moments and, in particular, for the variance. These results open the door for a discussion about the consistency of multiple tests, i.e.
	\begin{align}\label{eqn:def_consist}
		\frac{V_m}{R_m} - \E\Bigl( \frac{V_m}{R_m} \Bigr)  \xrightarrow{P_m}0 \text{ or equivalently }\text{Var}\Bigl( \frac{V_m}{R_m} \Bigr)\to 0.
	\end{align}
Specific results are discussed in \Cref{sec:consis,sec:estimators}. If $\liminf_{m\to\infty} \text{FDR}_m>0$ then we have the following necessary condition for consistency:
\begin{align}\label{eqn:nec_cond_cons_P(V>0)}
	\lim_{m\to\infty}P_m(V_m>0)=1.
\end{align} 
As already stated we can not expect consistency in general. In the following we discuss the BH test for two extreme case.
\begin{example}\label{exam:intro_extreme_cases}
	\begin{enumerate}[(a)]
		\item\label{enu:exam:intro_extreme_cases_DU} Let $m_1\geq 0$ be fixed. Then $P_m(V_m^{\text{BH}}=0)$ is minimal for the so-called Dirac uniform configuration DU$(m,m_1)$, where all entries of $(p_{i,m})_{i\in I_{1,m}}$ are equal to zero. Under this configuration  $V_m^{\text{BH}}(\alpha,m_1) \to V_{\text{SU}}(\alpha,m_1)$ in distribution with 
		\begin{align*}
			P\Bigl( V_{\text{SU}}(\alpha,m_1)=0 \Bigr)=(1-\alpha)\exp(-m_1\alpha),
		\end{align*} 
		see \citet{FinnerRoters2001} and Theorem 4.8 of \citet{Scheer2012}. The limit variable belongs to the class of linear Poisson distributions, see \citet{FinnerETAL2015}, \citet{Jain1975} and \citet{ConsulFamoye2006}. Hence, the BH-test is never consistent under BI for fixed $m_1$ since \eqref{eqn:nec_cond_cons_P(V>0)} is violated.
		
		\item\label{enu:exam:intro_extreme_cases_iid} Another extreme case is given by i.i.d. distributed $p$-value $(p_{i,m})_{i\in I_{m,1}}$. Suppose that $p_{i,m}$, $i\in I_{1,m}$, is uniformly distributed on $[0,\lambda]$, where $\alpha<\lambda<1$. Then the BH-tests are not consistent, see \Cref{theo:nec_for_cons}\eqref{enu:theo:nec_for_cons_unif_(0,lamb)}.
	\end{enumerate}
\end{example}
More information about DU$(m,m_1)$ and least favorable configurations can be found in \Cref{sec:LFC}. 
The requirement for consistency will be somehow in between these two extreme cases where the assumption $m_1\to\infty$ will be always needed. 

\subsection{Our step-up tests and underlying assumptions}
\begin{figure}[ht]
	\begin{center}
		\includegraphics[trim = 0mm 79.5mm 0mm 88mm,clip, scale =0.33]{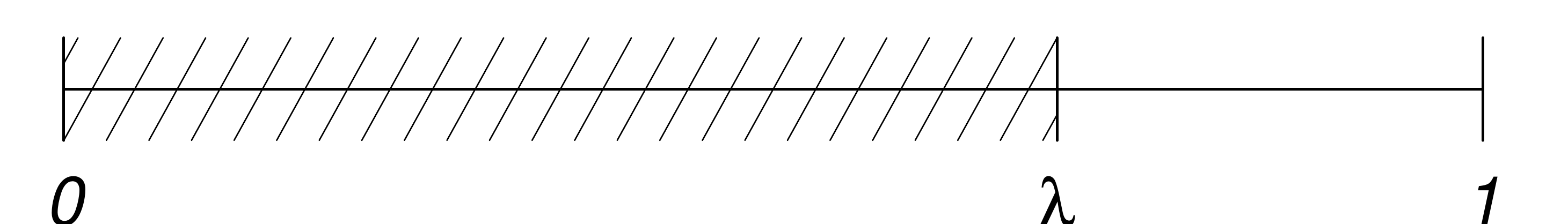}
	\end{center}
	\caption{Decision region $[0,\lambda]$ (dashed) and estimation region $(\lambda,1]$.}
	\label{fig:estimation_region}
\end{figure}

In the following we introduce the adaptive step-up tests we consider in this paper.  Let $0<\alpha<1$ be a fixed level and let $\lambda$, $\alpha\leq \lambda < 1$, be a tuning parameter and we agree that no null ${\mathcal H }_{i,m}$ with $p_{i,m}>\lambda$ should be rejected. The latter is not restrictive for practice since it is very unusual to reject a null if the corresponding $p$-value exceeds, for instance, $1/2$. We divide the range $[0,1]$ of the $p$-values in a decision region $[0,\lambda]$, where all ${\mathcal H }_{i,m}$ with $p_{i,m}\leq \lambda$ have a chance to be rejected, and an estimation region $(\lambda,1]$, where $p_{i,m}> \lambda$ are used to estimate $m_0$, see \Cref{fig:estimation_region}.
To be more specific we consider estimators $\widehat m_0$ of the form 
\begin{equation}\label{eqn:hatm0}
	\widehat m_0 = \widehat m_0((\widehat F_m(t))_{t\geq\lambda}) > 0
\end{equation}
for estimating $m_0$, which are measurable functions depending only on $(\widehat F_m(t))_{t\geq\lambda}$. As usual we denote by $\widehat F_m$ the empirical distribution function of the $p$-values $p_{1,m},\ldots,p_{m,m}$. As motivated in the introduction we now plug-in these estimators in the BH test. Doing this we obtain the data driven critical values
\begin{equation}\label{eqn:hat_alpha}
\widehat \alpha_{i:m} = \min\Bigl\{\left(\frac{i}{\widehat m_0} \alpha\right), \lambda\Bigr\}, \quad i=1,\ldots,n,
\end{equation}
where we promote to use the upper bound $\lambda$ as \citet{HeesenJanssen2016} already did. 
The following two quantities will be rapidly used:
Through
\begin{equation}\label{eqn:Rm+Vm}
	R_m(t) = m\widehat F_m(t) \quad \mbox{and} \quad V_m(t) = \sum_{i\in I_0} 1\{p_{i,m} \leq t\}, \;t\in[0,1].
\end{equation}
Throughout this paper, we investigate different mild assumptions. For our main results we fix the following two:
\begin{enumerate}[(\text{A}1)]
	\renewcommand\labelenumi{\bfseries(\text{A}\theenumi)}
	\item\label{enu:ass_A_kappa0} Suppose that 
	\begin{align*}
		\frac{m_0}{m}\to \kappa_0\in(0,1].
	\end{align*}
	
	\item\label{enu:ass_A_mo_leq_R} Suppose that $\widehat m_0$ is always positive and
	\begin{align*}
		\frac{\lambda}{\alpha}\widehat m_0 \geq R_m(\lambda).
	\end{align*}
	\xdef\tmp{\theenumi}
\end{enumerate}
If only  $0< \liminf_{m\to\infty}m_0/m$ is valid then our results apply to appropriate subsequences. The most interesting case is $\kappa_0>\alpha$ since otherwise (if $m_0/m\leq \alpha$) the FDR can be controlled, i.e. $\text{FDR}_m\leq \alpha$, by rejecting everything. 
\begin{remark}\label{rem:ass_A}
	\begin{enumerate}[(a)]
		\item %\label{enu:rem:ass_A_BI} 
		Under (A\ref{enu:ass_A_mo_leq_R}) the FDR of the adaptive multiple test was obtained for the BI model by \citet{HeesenJanssen2016}:
		\begin{align}\label{eqn:FDR_forumla}
			\text{FDR}_m=\frac{\alpha}{\lambda}\E \Bigl( \frac{V_m(\lambda)}{\widehat m_0} \Bigr).
		\end{align}
			In particular, we obtain
		\begin{align*}
			\text{FDR}_m \leq \E \Bigl( \frac{V_m(\lambda)}{R_m(\lambda)} \Bigr)\leq P_m(V_m(\lambda)>0),
		\end{align*}
		where the upper bound is always strictly smaller than $1$ for finite $m$.
		 
		\item\label{enu:rem:ass_A_doublehat} If (A\ref{enu:ass_A_mo_leq_R}) is not fulfilled then consider the estimator $\max\{\widehat m_0,(\alpha/\lambda)R_m(\lambda)\}$ instead of $\widehat m_0$. Note that both estimators lead to the same critical values $\alpha_{i:m}$ and so the assumption (A\ref{enu:ass_A_mo_leq_R}) is not restrictive.
	\end{enumerate}
\end{remark}
A prominent example for an adaptive test controlling the FDR by $\alpha$ is given by the Storey estimator \eqref{eqn:Storey_estimator_tilde_m0}: 
\begin{align}
	&{\widehat m}_0^{\text{Stor}}(\lambda):=\min\{\widetilde m_0(\lambda),\frac{\alpha}{\lambda}R_m(\lambda)\}\text{ with } \label{eqn:Storey_estimator}\\
	&\widetilde m_0^{\text{Stor}}(\lambda) = m \frac{ 1 - \widehat F_m(\lambda) + \frac{1}{m} }{1-\lambda}. \label{eqn:Storey_estimator_tilde_m0}
\end{align}
A refinement was established by \citet{HeesenJanssen2016}. They introduced a couple of inspection points $0<\lambda=\lambda_0<\lambda_1<\ldots< \lambda_k=1$, where $m_0$ is estimated on each interval $(\lambda_{i-1},\lambda_i]$. As motivation for this idea observe that the Storey estimator can be rewritten as the following linear combination
\begin{align*}
	\widetilde m_0^{\text{Stor}}(\lambda) = \sum_{i=1}^k \beta_i m \frac{ \widehat F_m(\lambda_i) - \widehat F_m(\lambda_{i-1}) + \frac{1}{m} }{\lambda_i-\lambda_{i-1}}
\end{align*}
with weights $\beta_i=(\lambda_i-\lambda_{i-1})/(1-\lambda)$, where $\sum_{i=1}^k\beta_i=1$. The ingredients 
\begin{align}\label{eqn:generl_Storey_estimator_tilde_m0}
	\widetilde m_0(\lambda_{i-1},\lambda_i) := m \frac{ \widehat F_m(\lambda_i) - \widehat F_m(\lambda_{i-1}) + \frac{1}{m} }{\lambda_i-\lambda_{i-1}} 
\end{align}
are also estimators for $m_0$, which were used by \citet{LiangNettleton2012} in another context. Under BI the following theorem was proved by \citet{HeesenJanssen2016}.  A discussion of their consistency is given in \Cref{sec:estimators}.
\begin{theorem}[cf. Thm 10 in \cite{HeesenJanssen2016}]\label{theo:Storey_comb_FDR}
	Let $\widehat\beta_{i,m}=\widehat\beta_{i,m}((\widehat F_m(t))_{t\geq \lambda_i})\geq 0$ be random weights for $i\leq k$ with $\sum_{i=1}^k\widehat\beta_{i,m}=1$. The adaptive step-up tests using the estimator
	\begin{align}\label{eqn:gen_storey}
		\widetilde m_0:=\sum_{i=1}^k \widehat \beta_{i,m} \widetilde m_0(\lambda_{i-1},\lambda_i)
	\end{align}
	controls the FDR, i.e. $\text{FDR}_m\leq \alpha$.                                                                                                                                                                                                            
\end{theorem}
Finally, we want to present a necessary condition of asymptotic FDR control. It was proven by \citet{HeesenJanssen2015} for a greater class than the BI models, namely reverse martingals. The same condition was already used by \citet{FinnerGontscharuk2009} for asymptotic FWER control.
\begin{theorem}[cf. Thm 6.1 in \cite{HeesenJanssen2015}]\label{theo:nesc_cond_asymp_control}
	Suppose that (A\ref{enu:ass_A_kappa0}), (A\ref{enu:ass_A_mo_leq_R}) holds. If
	\begin{align*}
	P_m\Bigl( \frac{\widehat m_0}{m_0}\leq 1-\delta \Bigr)\to 0\text{ for all }\delta>0
\end{align*}
then we have asymptotic FDR control, i.e. $\limsup_{m\to\infty}\text{FDR}_m\leq \alpha$.
\end{theorem}

\section{Moments}\label{sec:moments}
This section provides exact second moment formulas of $\text{FDP}_m=V_m/{R_m}$ for our adaptive step-up tests for a fixed regime $P_m$. Our method of proof relies on conditioning with respect to the $\sigma$-algebra
\begin{align*}
	\mathcal{F}_{\lambda,m} := \sigma\Bigl( \mathbf{1}\{p_{i,m}\leq s\}: s\in[\lambda,1], \, 1 \leq i \leq m \Bigr).
\end{align*}
Conditional under the (non-observable) $\sigma$-algebra $\mathcal{F}_{\lambda,m}$ the quantities $\widehat m_0$, $R_m(\lambda)$ and $V_m(\lambda)$ are fixed values. But only $R_m(\lambda)=m\widehat F_m(\lambda)$ and $\widehat m_0$ are given by the data and observable. The FDR formula \eqref{eqn:FDR_forumla} is now completed by an exact variance formula. The proof offers also a rapid approach to the known variance formula of \citet{FerreiraZwinderman2006} for the Benjamini and Hochberg test (with $\widehat m_0=m$ and $\lambda=\alpha$). Without loss of generality we can assume that $p_m=(p_{1,m},\ldots,p_{m,m})$ is ordered by
\begin{align*}
	I_{0,m}=\{1,\ldots,m_0\}\text{ and }I_{1,m}=\{m_0+1,\ldots,m\}.
\end{align*}
Now, we introduce a new $p$-value vector $p_m^{(1,\lambda)}$. If $V_m(\lambda)>0$ then set $p_m^{(1,\lambda)}$ equal to $p_m$ but replace one $p$-value $p_{i,m}$ with $p_{i,m}\leq \lambda$ by $0$ for one $i\leq m_0$, for convenience take the smallest integer $i$ with this property. If $V_m(\lambda)=0$ then set $p_m^{(1,\lambda)}=p_m$. Moreover, let $R_m^{(1,\lambda)}=R_m^{(1,\lambda)}(p_m^{(1,\lambda)})$ be the number of rejections of the adaptive test for substituted vector $p_m^{(1,\lambda)}$ of $p$-value. Note that $\widehat m_0$ remains unchanged when considering $p_m^{(1,\lambda)}$ instead of $p$.
\begin{theorem}\label{theo:form_moments}
	Suppose that our assumptions (A\ref{enu:ass_A_mo_leq_R}) are fulfilled: 
	\begin{enumerate}[(a)]
		\item \label{enu:theo:form_moments_second_moment} The second moment of $\text{FDP}_m$ is given by
		\begin{align*}
			\E \Bigl( \Bigl( \frac{V_m}{R_m} \Bigr)^2 \Bigr)= \E \Bigl( \frac{\alpha^2V_m(\lambda)(V_m(\lambda)-1)}{ \lambda^2\widehat m_0^2} + \frac{\alpha}{\lambda} \frac{ V_m(\lambda)}{\widehat m_0} \E \Bigl( \frac{1}{R_m^{(1,\lambda)}} \Bigl | \mathcal{F}_{\lambda,m} \Bigr) \Bigr).
		\end{align*}
		
		\item \label{enu:theo:form_moments_variance} The variance of $\text{FDP}$ fulfils
		\begin{align*}
			\text{Var}\Bigl( \frac{V_m}{R_m} \Bigr) = \frac{\alpha^2}{\lambda^2}\Bigl[\frac{\lambda}{\alpha} \E \Bigl( \frac{V_m(\lambda)}{\widehat m_0 }\E \Bigl( \frac{1}{R_m^{(1,\lambda)}} \Bigl | \mathcal{F}_{\lambda,m} \Bigr) \Bigr)+ \text{Var}\Bigl( \frac{V_m(\lambda)}{\widehat m_0} \Bigr) -  \E \Bigl( \frac{V_m(\lambda)}{\widehat m_0^2} \Bigr) \Bigr].
		\end{align*}
		
		\item \label{enu:theo:form_moments_E(V)} We have
		\begin{align*}
			 \E(V_m)= \frac{\alpha}{\lambda}\E \Bigl( \frac{V_m(\lambda)}{\widehat m_0} \E(R_m^{(1,\lambda)} |\mathcal{F}_{\lambda,m} )\Bigr).
		\end{align*}
	\end{enumerate}
\end{theorem}
Exact higher moment formulas are established in the appendix, see \Cref{sec:higher_moments}.

\section{The variability of $\text{FDP}_m$ and the Consistency of adaptive multiple tests}\label{sec:consis} 
The exact variance formula applies to the stability of the $\text{FDP}$ and its consistency if $m$ tends to infinity. If not stated otherwise, all limits are meant as $m\to\infty$. In the following we need a further mild assumption:
\begin{enumerate}[1]
	\setcounter{enumi}{\tmp} % Dadurch wird der Z{\"a}hler der letzten Aufz{\"a}hlung {\"u}bernommen
	\renewcommand\labelenumi{\bfseries(\text{A}\theenumi)}
	\item\label{enu:ass_A_mo_leq_Km} There is some $K>0$ such that $\widehat m_0\leq Km$ for all $m\in\N$.
	\xdef\tmp{\theenumi}
\end{enumerate}
Clearly, (A\ref{enu:ass_A_mo_leq_Km}) is fulfilled for the trivial estimator $\widehat m_0=m$ and for all generalized weighted estimators of the form \eqref{eqn:gen_storey} with $K=2\sum_{i=1}^k(\lambda_i-\lambda_{i-1})^{-1}$. Note that (A\ref{enu:ass_A_kappa0}) and (A\ref{enu:ass_A_mo_leq_Km}) imply $\liminf_{m\to\infty}\text{FDR}_m>0$ and, hence, \eqref{eqn:nec_cond_cons_P(V>0)} is a necessary condition for consistency in this case.  In the following we give boundaries for the variance of $\text{FDP}_m=V_m/R_m$ depending on the leading term in the variance formula of \Cref{theo:form_moments}:
\begin{align*}
C_{m,\lambda}:=\frac{\alpha}{\lambda}\E \Bigl( \frac{V_m(\lambda)}{\widehat m_0} \E \Bigl( \frac{1}{R_m^{(1,\lambda)}}\Bigl |{\mathcal F }_{\lambda,m} \Bigr) \Bigr)+\Bigl( \frac{\alpha}{\lambda} \Bigr)^2\text{Var}\Bigl( \frac{V_m(\lambda)}{\widehat m_0} \Bigr).
\end{align*}
\begin{lemma}\label{lem:bound_var}
	Suppose that (A\ref{enu:ass_A_mo_leq_R}) is fulfilled.
	\begin{enumerate}[(a)]
		\item \label{enu:lem_bound_var_Cm} We have
		\begin{align}
			&\E\Bigl( \frac{V_m(\lambda)}{\widehat m_0^2} \Bigr)\leq  \Bigl( \frac{\lambda}{\alpha} \Bigr)^2 \frac{2}{\lambda(m_0+1)}\text{ and }\label{eqn:E(V/m0^2)} \\
			&C_{m,\lambda} \leq \text{Var}\Bigl( \frac{V_m}{R_m} \Bigr)\leq C_{m,\lambda} + \frac{2}{\lambda(m_0+1)}.\label{eqn:bounds_variance}
		\end{align}
				
		\item \label{enu:lem_bound_var_E(ind_E(V|F))} Suppose (A\ref{enu:ass_A_mo_leq_Km}). Then $\widehat m_0\leq m_0K_m$ with $K_m:=Km/m_0$ and for all $t>0$
		\begin{align*}%\label{eqn:bound_ind_E(V|F)}
			&P_m\Bigl( \E(V_m | \mathcal{F}_{\lambda,m} )\leq t\Bigr) \leq P_m(V_m=0) +tD_{m,\lambda},\text{ where}\\
			&D_{m,\lambda}:= \frac{ 4K_m^2}{\alpha^2}  \Bigl[ \text{Var}\Bigl( \frac{V_m}{R_m} \Bigr)+\frac{2}{\lambda(m_0+1)}-\frac{\alpha^2}{\lambda^2}\text{Var}\Bigl( \frac{V_m(\lambda)}{\widehat m_0} \Bigr)  \Bigr]+\frac{\lambda m_0 K_m}{\alpha \exp(\frac{1}{8}m_0\lambda)}.
		\end{align*}
	\end{enumerate}
\end{lemma}
Since under (A\ref{enu:ass_A_kappa0}) $m_0\to \infty$ we have consistency iff $C_{m,\lambda}\to 0$. In the following we present sufficient and necessary conditions for this.
\begin{theorem}\label{theo:consistency}
	Under (A\ref{enu:ass_A_kappa0})-(A\ref{enu:ass_A_mo_leq_Km}) the following \eqref{enu:theo:consis_consis} and \eqref{enu:theo:consis_var_R1} are equivalent.
	\begin{enumerate}[(a)]
		\item \label{enu:theo:consis_consis} (Consistency) We have $V_m/R_m-\E(V_m/R_m)\to 0$ in $P_m$-probability.
		
		\item \label{enu:theo:consis_var_R1} It holds that
		\begin{align}
			&\frac{\widehat m_0}{m}-\E\Bigl( \frac{\widehat m_0}{m_0} \Bigr)\to 0\text{ in }P_m\text{-probability} \text{ and } \label{eqn:theo:consistency_m0/m-E()} \\
			&R^{(1,\lambda)}_m\to \infty \text{ in }P_m\text{-probability}.\label{eqn:theo:consistency_R1}
		\end{align}
	\end{enumerate} 
\end{theorem}
Roughly speaking the consistency requires an amount of rejections \eqref{eqn:theo:consistency_R1} turning to infinity and a stability condition \eqref{eqn:theo:consistency_m0/m-E()} for the estimator $\widehat m_0$, which is equivalent to $\text{Var} (V_m(\lambda)/\widehat m_0 )\to 0$. 
\begin{remark}\label{rem:bound_var}
	Suppose that (A\ref{enu:ass_A_kappa0})-(A\ref{enu:ass_A_mo_leq_Km}) are fulfilled. 
\begin{enumerate}[(a)]
	\item $\text{Var}(V_m/R_m)\to 0$ implies 
	\begin{align}\label{eqn:rem_bound_var_varto0}
	\E \Bigl( \frac{1}{R_m^{(1,\lambda)}} \Bigl | \mathcal{F}_{\lambda,m} \Bigr) \to 0 \text{ in }P_m\text{-probability}.
	\end{align}
	
	\item\label{enu:rem:bound_var_E(V)_to_infinity} Under \eqref{eqn:rem_bound_var_varto0} we have $\E(R^{(1,\lambda)}_m | \mathcal{F}_{\lambda,m})\to\infty$ and $\E(V_m | \mathcal{F}_{\lambda,m})\to \infty$ in $P_m$-probability and so $\E(V_m)\to\infty$ by \Cref{theo:form_moments}\eqref{enu:theo:form_moments_E(V)}.	
\end{enumerate}
\end{remark}
Under mild additional assumptions we can improve the convergence in expectation $\E(V_m)\to \infty$ from \Cref{rem:bound_var}\eqref{enu:rem:bound_var_E(V)_to_infinity}. Recall that $V_m$ and $R_m$ depend, of course, on the pre-specified level $\alpha$. In comparison to the rest of this paper we consider in the following theorem we consider more than one level. That is why we prefer (only) for this theorem the notation $V_{m,\alpha}$ and $R_{m,\alpha}$
\begin{theorem}\label{theo:Vm_to_infty}                        
	Suppose (A\ref{enu:ass_A_kappa0})-(A\ref{enu:ass_A_mo_leq_Km}). Moreover, we assume that we have consistency for all level $\alpha\in(\alpha_1,\alpha_2)$ and some $0<\alpha_1<\alpha_2<1$. Then we have in $P_m$-probability for all $\alpha\in(\alpha_1,\alpha_2)$ that
	\begin{align*}
	V_{m,\alpha} \to \infty \text{ and so }R_{m,\alpha}\to \infty. 
	\end{align*}	
\end{theorem} 
The next example points out that consistency may depend on the level $\alpha$ and adaptive tests may be consistent while the BH test is not so. A proof of the statements is given in \Cref{sec:proof}.
\begin{example}\label{exam:consistency}
	Let $U_1,U_2,\ldots,U_m$ be i.i.d. uniformly distributed on $(0,1)$. Consider $1/2<\lambda<1$, $m_0=m_1$ and $p$-values from the false null given by $p_{i,m}=\min\{U_i, x_0\}$ with $x_0:=1/6$, $i\in I_{m,1}$.  The BH test BH($\alpha$) with level $\alpha:=1/4$ is not consistent while BH$(2\alpha)$ is consistent. But the adaptive test Stor$(\alpha,\lambda)$ using the classical Storey estimator \eqref{eqn:Storey_estimator} is consistent.
\end{example}

\section{Consistent and inconsistent regimes}\label{sec:estimators}
Below we will exclude the ugly estimator $\widehat m_0=(\alpha/\lambda)(R_m(\lambda)\vee 1)$ which could lead to rejecting all hypotheses with $p_{i,m}\leq \lambda$. To avoid this effect let us introduce:
\begin{enumerate}[(\text{A}1)]
	\setcounter{enumi}{\tmp} % Dadurch wird der Z{\"a}hler der letzten Aufz{\"a}hlung {\"u}bernommen
	\renewcommand\labelenumi{\bfseries(\text{A}\theenumi)}
	\item\label{enu:ass_A_mo_geq_CR} There exists a constant $C>1$ with
	\begin{align*}
	\lim_{m\to\infty}P_m\Bigl( \widehat m_0\geq \frac{C\alpha}{\lambda}(R_m(\lambda)\vee 1) \Bigr)=1.
	\end{align*}
	\xdef\tmp{\theenumi}
\end{enumerate}
Note that (A\ref{enu:ass_A_mo_geq_CR}) guarantees that (A\ref{enu:ass_A_mo_leq_R}) holds at least with probability tending to one. The next theorem yields a necessary condition for consistency.
\begin{theorem}\label{theo:nec_for_cons} Suppose that $\liminf_{m\to\infty}\text{FDR}_m>0$ and  (A\ref{enu:ass_A_mo_geq_CR}) holds.
	\begin{enumerate}[(a)]
		\item\label{enu:theo:nec_for_cons_m1_to_infty} If we have consistency then $m_1\to\infty$.
		
		\item\label{enu:theo:nec_for_cons_unif_(0,lamb)} Suppose that (A\ref{enu:ass_A_kappa0}) holds. If all $(p_{i,m})_{i\in I_{1,m}}$ are i.i.d. uniformly distributed on $[0,\lambda]$ then we have no consistency. 
	\end{enumerate}
\end{theorem}
Consistency for the case $\kappa_0<1$ was already discussed by \citet{GenoveseWassermann2004}, who used a stochastic process approach. Also \citet{FerreiraZwinderman2006} used their formulas for the moments of $FDP_m$ to discuss the consistency for the BH test. By their Proposition 2.2 or our \Cref{theo:consistency} it is sufficient to show for $\widehat m_0=m$ that $R_m^{BH}\to\infty$ in $P_m$-probability. For this purpose \citet{FerreiraZwinderman2006} found conditions such that $R_m/m\to \widetilde C>0$ in $P_m$-probability.  The sparse signal case $\kappa_0=1$ is more delicate since $R_m/m$ always tends to $0$ even for adaptive tests. Recall for the following lemma that $\widehat \alpha_{R_m:m}$ is the largest intersection point of $\widehat F_m$ and the random Simes line $t\mapsto f(t)= :(\widehat m_0/m)(t/\alpha)$, observe $\alpha_{i:m}=f^{-1}(i/m)$. 
\begin{lemma}\label{lem:R_m/m_to_0}
	Suppose that (A\ref{enu:ass_A_kappa0}) with $\kappa_0=1$ and (A\ref{enu:ass_A_mo_geq_CR}) are fulfilled. Then $\widehat \alpha_{R_m:m} \to 0$ in $P_m$-probability. In particular, under (A\ref{enu:ass_A_mo_leq_Km}) we have  $R_m/m\to 0$ in $P_m$-probability.
\end{lemma}
Besides the result of \Cref{theo:nec_for_cons} we already know that a further necessary condition for consistency is \eqref{eqn:theo:consistency_m0/m-E()} which is assumed to be fulfilled in the following. Turning to convergent subsequence we can assume without loss of generality under (A\ref{enu:ass_A_kappa0}), (A\ref{enu:ass_A_mo_leq_Km}) and (A\ref{enu:ass_A_mo_geq_CR}) that $\E(\widehat m_0/m)\to C_0\in[\kappa_0\alpha,K]$. In this case  \eqref{eqn:theo:consistency_m0/m-E()} is equivalent to
\begin{align}\label{eqn:theo:consistency_m0/m}
	\frac{\widehat m_0}{m}\to C_0\in[\kappa_0\alpha,K] \text{ in }P_m\text{-probability}. 
\end{align}
In the following we will rather work with \eqref{eqn:theo:consistency_m0/m} instead of  \eqref{eqn:theo:consistency_m0/m-E()}. Due to \Cref{lem:R_m/m_to_0} the question about consistency can be reduced in the sparse signal case $\kappa_0=1$ to the comparison of the random Simes line $f$ defined above and $\widehat F_m$ close to $0$.
\begin{theorem}\label{theo:cons_kappa=1}                                                                 
	Assume that (A\ref{enu:ass_A_kappa0}), (A\ref{enu:ass_A_mo_leq_Km}), (A\ref{enu:ass_A_mo_geq_CR}) and \eqref{eqn:theo:consistency_m0/m} hold. Let $\delta>0$ and $(t_m)_{m\in\N}$ be some sequence in $(0,\lambda)$  such that $mt_m\to \infty$ and
	\begin{align}\label{eqn:cons_kappa=1_condition}                                                                                   	P_m\Bigl( \frac{m_1}{m}\frac{ \widehat F_{1,m}(t_m) }{t_m} \geq \delta - \kappa_0 + \frac{1}{\alpha}C_0  \Bigr) \to 1,
	\end{align}
	where $\widehat F_{m,j}(x):=m_j^{-1}\sum_{i\in I_{m,j}}\mathbf{1}\{p_{i,m}\leq x\}, \,x\in(0,1)$, denotes for $j=0,1$ the empirical distribution function of the $p$-values corresponding to the true and false null, respectively. \\
	Then $V_m \to \infty$ in $P_m$-probability and so we have consistency by \Cref{theo:consistency}. 
\end{theorem}
\begin{remark}\label{rem:cons_kappa=1}
	\begin{enumerate}[(a)]
		\item \label{enu:rem:cons_kappa=1_F1}	Suppose that $(p_{i,m})_{i\in I_{1,m}}$ are i.i.d. with distribution function $F_1$. Then the statement of \Cref{theo:cons_kappa=1} remains valid if we replace \eqref{eqn:cons_kappa=1_condition} by the condition that for all sufficiently large $m\in\N$
		\begin{align}\label{eqn:cons_kappa=1_cond_F1}
			\frac{m_1}{m}\frac{ F_{1}(t_m) }{t_m} \geq \delta - \kappa_0 + \frac{1}{\alpha}C_0.
		\end{align}
	 	A proof of this statement is given in \Cref{sec:proof}. 
	 	
		\item In the case $\kappa_0=1$ we need a sequence $(t_m)_{m\in\N}$ tending to $0$. 
		
		\item For the   DU$(m_1,m)$-configuration the assumption \eqref{eqn:cons_kappa=1_condition}  is fulfilled for $t_m=(m_1/m)(K+2)$ as long as the necessary condition $m_1\to\infty$ holds. 
	\end{enumerate}
\end{remark}
As already stated, consistency only holds under certain additional assumptions. In the following we compare consistency of the classical BH test and adaptive tests with an appropriate estimator.
\begin{lemma}\label{lem:comp_BH_adapt_consis}
	Suppose that (A\ref{enu:ass_A_kappa0}), (A\ref{enu:ass_A_mo_leq_Km}) and (A\ref{enu:ass_A_mo_geq_CR}) are fulfilled. Assume that \eqref{eqn:theo:consistency_m0/m} holds for some $C_0\in[\alpha\kappa_0,1]$. If $C_0=1$ then additionally suppose that
	\begin{align}\label{eqn:lem_comp_BH_condition_to_1}
		P_m\Bigl( \frac{\widehat m_0}{m} \leq 1\Bigr)\to 1.
	\end{align}
	Then consistency of the BH test implies consistency of the adaptive test.
\end{lemma}
Under some mild assumptions \Cref{lem:comp_BH_adapt_consis} is applicable for the weighted estimator \eqref{eqn:gen_storey}, see \Cref{cor:Storey}\eqref{enu:cor:Storey_consis} for sufficient conditions.
\subsection{Combination of generalized Storey estimators}

In the following we become more concrete by discussing the combined Storey estimators \eqref{eqn:gen_storey} introduced in \Cref{sec:intro}. For this purpose we need the following assumption to ensure that (A\ref{enu:ass_A_mo_geq_CR}) is fulfilled. 
\begin{enumerate}[(\text{A}1)]                                                                                    	\setcounter{enumi}{\tmp} % Dadurch wird der Z{\"a}hler der letzten Aufz{\"a}hlung {\"u}bernommen
\renewcommand\labelenumi{\bfseries(\text{A}\theenumi)}
\item\label{enu:ass_comb_Storey_mo_geq_al_lam_R} Suppose that $\kappa_0 > \alpha(1-\kappa_0)/[\lambda(1-\alpha)]$.
\xdef\tmp{\theenumi}
\end{enumerate}

\begin{corollary}\label{cor:Storey}
Let (A\ref{enu:ass_A_kappa0}), (A\ref{enu:ass_comb_Storey_mo_geq_al_lam_R}) and the assumptions of \Cref{theo:Storey_comb_FDR} be fulfilled. Consider the  adaptive multiple test with $\widehat m_0=\widetilde m_0\vee (\alpha/\lambda)R_m(\lambda)$.
\begin{enumerate}[(a)]
	\item\label{enu:cor:Storey_ka=1} Suppose that $\kappa_0=1$. Then \eqref{eqn:theo:consistency_m0/m} holds with $C_0=1$ and
	\begin{align}\label{eqn:cor:Storey_ka=1_FDR}
	\lim_{m\to\infty}\text{FDR}_m= \alpha = \lim_{m\to\infty} \text{FDR}_m^{\text{BH}}.
	\end{align} 
	
	\item \label{enu:cor:Storey_ka<1} Suppose that 	$\kappa_0<1$ and we have with probability one that
	\begin{align}\label{eqn:cor:Storey_betai}
	(\lambda_1-\lambda_{0})^{-1}\widehat  \beta_{1,m}\leq(\lambda_2-\lambda_1)^{-1}\widehat \beta_{2,m} \leq \ldots\leq (1-\lambda_{k-1})^{-1}\widehat \beta_{m,m}
	\end{align}
	for every $m\in\N$. Moreover, assume that 
	\begin{align}\label{eqn:cor:Storey_condition_liminf}
		\liminf_{m\to\infty}\widehat F_{m}(\lambda_i)\geq\lambda_i + \varepsilon_{i}\text{ a.s. for some }\varepsilon_i\in[0,1-\lambda_i]
	\end{align}
	and all $i=1,\ldots,k$. If there is some $j\in\{1,\ldots,k\}$ and $\delta>0$ such that 
	\begin{align}\label{eqn:cor:Storey_condition_eps>0}
		\liminf_{m\to\infty}\frac{\widehat \beta_{j,m}}{\lambda_{j}-\lambda_{j-1}}- \frac{\widehat\beta_{j-1,m} }{\lambda_{j-1}-\lambda_{j-2}}\geq \delta \text{ a.s. and }\varepsilon_{j}>0,
	\end{align}
	where $\widehat \beta_0:=0=:\lambda_{-1}$, then we have an improvement of the $\text{FDR}_m$ asymptotically compared to the Benjamini-Hochberg procedure, i.e.
	\begin{align}\label{eqn:cor:Storey:improvement}
	\liminf_{m\to\infty}\text{FDR}_m> \kappa_0 \alpha = \lim_{m\to\infty} \text{FDR}_m^{\text{BH}}.
	\end{align}
	
	\item\label{enu:cor:Storey_consis} (Consistency) Suppose that the weights are asymptotically constant, i.e. $\widehat \beta_{i,m}\to\beta_{i}$ a.s. for all $i\in\{1,\ldots,k\}$, and fulfil \eqref{eqn:cor:Storey_betai}. Assume that 
	\begin{align}\label{eqn:cor:Storey_condition_lim}
		\lim_{m\to\infty}\widehat F_{m,1}(\lambda_i)=\lambda_i + \varepsilon_{i}\text{ a.s. for some }\varepsilon_i\in[0,1-\lambda_i]
	\end{align}
	and for all $i=1,\ldots,k$. Moreover, suppose that 
	\begin{align}\label{eqn:cor:Storey_condition_eps>0_consis}
		\gamma_j:=\frac{\beta_{j}}{\lambda_{j}-\lambda_{j-1}}- \frac{\beta_{j-1} }{\lambda_{j-1}-\lambda_{j-2}}>0\,\text{ and }\,\varepsilon_{j}>0
	\end{align}
	for some $j\in\{1,\ldots,k\}$, 	where $\beta_0:=\widehat \beta_0:=0=:\lambda_{-1}$. Additionally, assume $m_1/\sqrt{m}\to\infty$ if $\kappa_0=1$. Then \eqref{eqn:theo:consistency_m0/m} holds for some $C_0\in[0,1]$ and consistency of the BH test implies consistency of the adaptive test. Moreover, if \eqref{eqn:cons_kappa=1_condition} holds for $C_0$ and a sequence $(t_m)_{m\in\N}$ with $mt_m\to\infty$ then we have always consistency of the adaptive test.
\end{enumerate}
\end{corollary}                                                                                 
It is easy to see that the assumptions of \eqref{enu:cor:Storey_consis} imply the ones of \eqref{enu:cor:Storey_ka<1}.                                                                                                                   
Typically, the $p$-values $p_{i,m}$, $i\in I_{1,m}$, from the false null are stochastically smaller than the uniform distribution, i.e. $P_m(p_{i,m}\leq x) \geq x$ for all $x\in(0,1)$ (with strict inequality for some $x=\lambda_i$). This may lead to \eqref{eqn:cor:Storey_condition_liminf} or \eqref{eqn:cor:Storey_condition_lim}.
\begin{remark}\label{rem:cor:Storey}
	If $p_{m_0+1,m},\ldots,p_{m,m}$ are i.i.d. with distribution function $F_1$ such that $F_1(\lambda_i)\geq\lambda_i$ for all $i=1,\ldots,k$. Then \eqref{eqn:cor:Storey_condition_liminf} and \eqref{eqn:cor:Storey_condition_lim} are fulfilled. Moreover, if $\kappa_0<1$ and $F_1(\lambda_i)>\lambda_i$ then $\varepsilon_{i}>0$. 
\end{remark}
If the weights $\widehat\beta_i=\beta_i$ are deterministic then weights fulfilling \eqref{eqn:cor:Storey_betai} produce convex combinations of Storey estimators with different tuning parameters $\lambda_i$, compare to \eqref{eqn:Storey_estimator}-\eqref{eqn:generl_Storey_estimator_tilde_m0}. 

\subsection{Asymptotically optimal rejection curve}\label{sec:AORC}
Our results can be transferred to general deterministic critical values \eqref{stepup001}, which are not of the form \eqref{eqn:hat_alpha} and do not use a plug-in estimator for $m_0$. To label this case we use $\lambda=1$. Analogously to \Cref{sec:consis} and \Cref{sec:higher_moments} we define $R_m^{(j,1)}$ for $j\in\N$ by setting $j$ $p$-values from the true null to $0$. By the same arguments as in the proof of \Cref{theo:form_moments} we obtain
\begin{align*}
	&E\Bigl( \frac{V_m}{R_m} \Bigr) = m_0 \E \Bigl( \,\frac{\alpha_{R_m^{(1,1)}:m}}{R_m^{(1,1)}}\, \Bigr)\text{ and }\\
	&E\Bigl( \Bigl( \frac{V_m}{R_m} \Bigr)^2 \Bigr)=m_0  \E \Bigl( \,\frac{\alpha_{R_m^{(1,1)}:m}}{ (R_m^{(1,1)})^2} \,\Bigr)+m_0(m_0-1)  \E \Bigl(\Bigl(\, \frac{ \alpha_{R_m^{(2,1)}:m}}{ (R_m^{(2,1)})}\, \Bigr)^2  \Bigr).
\end{align*}
The first formula can also be found in \citet{BenditkisETAL2018}, see the proof of Theorem 2 therein. The proof of the second one is left to the reader. By these formulas we can now treat an important class of critical values given by
\begin{align}\label{eqn:alpha_im_quotient}
	\alpha_{i:m}=\frac{i\alpha}{m+b-ai},\, i\leq m, \, 0\leq \min(a,b).
\end{align}
A necessary condition for the valid step-up tests is $\alpha_{m:m}<1$. This condition holds for the critical values \eqref{eqn:alpha_im_quotient} if
\begin{align}\label{eqn:nec_cond_alpha_quot}
	b>0\text{ and }a\in[0, 1-\alpha] \text{ or }b=0 \text{ and }a\in[0, 1-\alpha).
\end{align}
These critical values are closely related to
\begin{align}\label{eqn:alpha_AORC}
	\alpha^{\text{AORC}}_{i:m}=\frac{i\alpha}{m-i(1-\alpha)}=f_\alpha^{-1}\Bigl( \frac{i}{m} \Bigr),\,i<m,
\end{align}
of the asymptotically optimal rejection curve $f_\alpha(t)=t/(t(1-\alpha)+\alpha)$ introduced by \citet{FinnerETAL2009}. Note that the case $i=m$ is excluded on purpose because it would lead to $\alpha^{\text{AORC}}_{m:m}=1$. The remaining coefficient $\alpha^{\text{AORC}}_{m:m}$ has to be defined separately such that $\alpha^{\text{AORC}}_{m-1:m}\leq \alpha^{\text{AORC}}_{m:m}<1$, see \citet{FinnerETAL2009} and \citet{Gontscharuk2010} for a detailed discussion. It is well-known that neither for \eqref{eqn:alpha_im_quotient} with $b=0$ and $a>0$ nor for \eqref{eqn:alpha_AORC} we have control of the FDR by $\alpha$ over all BI models simultaneously. This follows from Lemma 4.1 of \citet{HeesenJanssen2015} since $\alpha_{1:m}>\alpha/m$. However, \citet{HeesenJanssen2015} proved that for all fixed $b>0$, $\alpha\in(0,1)$ and $m\in\N$ there exists a unique parameter $a_m\in(0,b)$ such that
\begin{align*}
	\sup_{P_m} \text{FDR}_{(b,a_m)}=\alpha,
\end{align*}
where the supremum is taken over all BI models at sample size $m$. The value $a_m$ may be found under the least favorable configuration DU$(m,m_1)$ using numerical methods.\\
By transferring our techniques to this type of critical values we get the following sufficient and necessary conditions for consistency.
\begin{lemma}\label{lem:cons_gen_alp}
	Let (A\ref{enu:ass_A_kappa0}) be fulfilled. Let $(a_m)_{m\in\N}$ and $(b_m)_{m\in\N}$ be sequences in $\R$ such that $b_m/ m\to 0$ and $(a_m,b_m)$ fulfil \eqref{eqn:nec_cond_alpha_quot} for every $m\in\N$. For every $m\in\N$ consider the step-up test with critical values given by \eqref{eqn:alpha_im_quotient} with $(a,b)=(a_m,b_m)$.
	\begin{enumerate}[(a)]
		\item \label{enu:lem:cons_gen_alp_suff_cond} Then we have consistency, i.e. $V_m/R_m-E(V_m/R_m)\to 0$ in $P_m$-probability, iff the following conditions \eqref{eqn:lem_cons_gen_alp_R1_to_infty}-\eqref{eqn:lem:cons_gen_alp_def_psi1} hold in $P_m$-probability: 
		\begin{align}
			&R_m^{(1,1)}\to \infty,\label{eqn:lem_cons_gen_alp_R1_to_infty}\\
			&\frac{a_m}{m}\Bigr(R_m^{(2,1)}-R_m^{(1,1)}\Bigl)\to 0 ,\label{eqn:lem_cons_gen_alp_R2-R1} \\
			&\frac{m+b_m-a_mR_m^{(1,1)}}{m_0} - \E\Bigl( \frac{m+b_m-a_mR_m^{(1,1)}}{m_0} \Bigr)\to 0. \label{eqn:lem:cons_gen_alp_def_psi1}
		\end{align}
		
		\item \label{enu:lem:cons_gen_alp_suff_cond_sparse} If $\kappa_0=1$, $m_1\to\infty$ and $\limsup_{m\to\infty}a_m<1-\alpha$ then \eqref{eqn:lem_cons_gen_alp_R1_to_infty} is sufficient for consistency and, moreover, $E(V_m/R_m)\to \alpha$ in this case.
	\end{enumerate}
\end{lemma}

\section{Least favorable configurations and consistency}\label{sec:LFC}
Below least favorable configurations (LFC) are derived for the $p$-value $(p_{i,m})_{i\in I_{m,1}}$ of the false portion. When deterministic critical values $i\mapsto \alpha_{i:m}/i$ are increasing then the FDR is decreasing in each argument $p_{i,m}$, $i\in I_{1,m}$, for fixed $m_1$, see \citet{BenjaminiYekutieli2001} or \citet{BenditkisETAL2018} for a short proof. Here and subsequently, we use "increasing" and "decreasing" in their weak form, i.e. equality is allowed, whereas other authors use "nondecreasing" and "nonincreasing" for this purpose. In that case the Dirac uniform configuration DU$(m,m_1)$, see \Cref{exam:intro_extreme_cases}, has maximum FDR, i.e. it is LFC. LFC are sometimes useful tools for all kind of proofs.

\begin{remark}\label{rem:LFC}
	In contrast to \eqref{eqn:hat_alpha} the original Storey adaptive test is based on $\widehat \alpha^{\text{Stor}}_{i:m}=(i/\widetilde m_0)\alpha$ for the estimator from \eqref{eqn:generl_Storey_estimator_tilde_m0}. It is known that in this situation DU$(m,m_1)$ is not LFC for the FDR, see \citet{BlanchardETAL2014}. However, we will see that for our modification $\widehat \alpha^{\text{Stor}}_{i:m}\wedge \lambda$ the DU$(m,m_1)$-model is LFC.
\end{remark}
Our exact moment formulas provide various LFC-results which are collected below. To formulate these we introduce a new assumption
\begin{enumerate}[(\text{A}1)]
	\setcounter{enumi}{\tmp} % Dadurch wird der Z{\"a}hler der letzten Aufz{\"a}hlung {\"u}bernommen
	\renewcommand\labelenumi{\bfseries(\text{A}\theenumi)}
	\item\label{enu:ass_A_mo_increasing} Let $p_{j,m}\mapsto \widehat m_0((p_{i,m})_i\leq m)$ be increasing for each coordinate $j\leq m$. 
	\xdef\tmp{\theenumi}
\end{enumerate}
Below we are going to condition on $(p_{i,m})_{i\in I_{1,m}}$. By (BI\ref{BI2}) we may write $P_m=P_{0,m}\otimes P_{1,m}$, where $P_{j,m}$ 
represents the distribution of $(p_{i,m})_{i\in I_{j,m}}$ under $P_m$ for $j\in\{0,1\}$, and $E(X|((p_{i,m})_{i\in I_{1,m}}))=\int X((p_{i,m})_{i\leq m})\,\mathrm{ d }P_{0,m}((p_{i,m})_{i\in I_{0,m}})$.
\begin{theorem}[LFC for adaptive tests]\label{theo:LFC_ada_test}
	Suppose that (A\ref{enu:ass_A_mo_leq_R}) is fulfilled. Define the vector   $p^*_{\lambda,m}:=(p_{i,m}\mathbf{1}\{p_{i,m}>\lambda\})_{i\in I_{1,m}}$. 
	\begin{enumerate}[(a)]
		\item\label{enu:theo:LFS_ada_test_cond_LFC} (Conditional LFC) 
		\begin{enumerate}[(i)]
			\item\label{enu:theo:LFS_ada_test_cond_LFC_Exp} The conditional FDR conditioned on $(p_{i,m})_{i\in I_{1,m}}$
			\begin{align*}
				\E \Bigl( \frac{V_m}{R_m} \Bigl |(p_{i,m})_{i\in I_{1,m}}\Bigr)=\E \Bigl( \frac{V_m}{R_m} \Bigl |p^*_{\lambda,m}\Bigr)
			\end{align*}
			only depends on the portion $p_{i,m}>\lambda$, $i\in I_{1,m}$.
			
			\item\label{enu:theo:LFS_ada_test_cond_LFC_var} Conditioned on $p^*_{\lambda,m}$ a configuration $(p_{i,m})_{i\in I_{1,m}}$ is conditionally Dirac uniform if $p_{i,m}=0$ for all $p_{i,m}\leq \lambda$, $i\in I_{1,m}$. The conditional variance of $V_m/R_m$ 
			\begin{align*}
				\text{Var}\Bigl( \frac{V_m}{R_m} \Bigl | p^*_{\lambda,m} \Bigr):= \E \Bigl( \Bigl( \frac{V_m}{R_m}\Bigr)^2 \Bigl | p^*_{\lambda,m}\Bigr)-\E \Bigl( \frac{V_m}{R_m}\Bigl | p^*_{\lambda,m}\Bigr)^2 
			\end{align*}
			is minimal under DU$_{\text{cond}}(m,M_{1,m}(\lambda))$, where $M_{1,m}(\lambda):=R_m(\lambda)-V_m(\lambda)$ is fixed conditionally on $p^*_{\lambda,m}$.
			
		\end{enumerate}
	\item\label{enu:theo:LFS_ada_test_uncond_LFC} (Comparison of different regimes $P_{1,m}$) Under  (A\ref{enu:ass_A_mo_increasing}) we have:
	\begin{enumerate}[(i)]
		\item\label{enu:theo:LFS_ada_test_uncond_LFC_decrea} If $p_{i,m}$ decreases for some $i\in I_{1,m}$ then $\text{FDR}_m$ increases. \\
		If $p_{i,m}\leq \lambda$,  $i\in I_{1,m}$, decreases then Var$_m(\text{FDP}_m)$ decreases. 
		
		\item\label{enu:theo:LFS_ada_test_uncond_LFC_max+min_DU} For fixed $m_1$ the DU$(m,m_1)$ configuration is LFC with maximal $\text{FDR}_m$. Moreover, it has minimal Var$_m(\text{FDP}_m)$ for all models with $p_{i,m}\leq \lambda$ a.s. for all $i\in I_{1,m}$. 
	\end{enumerate}
	\end{enumerate}
\end{theorem}
While any deterministic convex combination of Storey estimators $\widetilde m^{\text{Stor}}_0(\lambda_i)$ fulfils (A\ref{enu:ass_A_mo_increasing}) it may fail for estimators of the form \eqref{eqn:generl_Storey_estimator_tilde_m0}. But if the weights fulfil \eqref{eqn:cor:Storey_betai} then (A\ref{enu:ass_A_mo_increasing}) holds also for a convex combination \eqref{eqn:gen_storey} of these estimators. This follows from the other representation of the estimator \eqref{eqn:gen_storey} used in the proof of \Cref{cor:Storey}\eqref{enu:cor:Storey_consis}.

\section{Discussion and summary}\label{sec:discussion}
In this paper we presented finite sample variance and higher moments formulas for the false discovery proportion (FDP) of adaptive step-up tests. These formulas allow a better understanding of FDP. Among others, the formulas can be used to discuss consistency of FDP, which is preferable for application since the fluctuation and so the uncertainty vanishes. We determined a sufficient and necessary two-part condition for consistency: 
\begin{enumerate}[(i)]
	\item We need a stable estimator in the sense that $\widehat m_0/m-E(\widehat m_0/m)$ tends to 0 in probability.
	
	\item\label{enu:discu_false_p}  The $p$-values of the false null need to be stochastically small "enough" compared to the uniform distribution such that the number of rejections tends to $\infty$ in probability.
\end{enumerate}
Since the latter is more difficult to verify we gave a sufficient condition for it, see \eqref{eqn:cons_kappa=1_condition}. This condition also applies to the sparse signal case $m_0/m\to 0$, which is more delicate than the usual studied case $m_0/m\to K>0$. \\
In addition to the general results we discussed data dependently weighted combinations of generalized Storey estimators. Tests based on these estimators were already discussed by \citet{HeesenJanssen2016}, who showed finite FDR control by $\alpha$. In \citet{Heesen2014} and \citet{HeesenJanssen2016} there are practical guidelines how to choose the data dependent weights. But note that for our results the weights have to fulfil the additional condition \eqref{eqn:cor:Storey_betai}. We want to summarize briefly advantages of these tests in comparison to the classical BH test (compare to \Cref{cor:Storey}\eqref{enu:cor:Storey_consis}):
\begin{itemize}
	\item The adaptive tests attain (if $\kappa_0=1$) or even exhaust (if $\kappa_0<1$) the (asymptotic) FDR level $\kappa_0\alpha$ of the BH test.
	
	\item Under mild assumptions consistency of the BH test always implies consistency of the adaptive test.
\end{itemize}
In \Cref{sec:AORC} we explained that our results can also be transferred to general deterministic critical values $\alpha_{i:m}$, which are not based on plug-in estimators of $m_0$. The same should be possible for general random critical values under appropriate conditions. Due to lack of space we leave a discussion about other estimators for future research.

\section{Proofs}\label{sec:proof}
\subsection{Proof of \Cref{theo:form_moments}}
To improve the readability of the proof, all indices $m$ are submitted, i.e. we write $p_i$ instead of $p_{i,m}$ etc. First, we determine $\E(\text{FDP}^2|{\mathcal F }_{\lambda})$. Without loss of generality we can assume conditioned on $\mathcal{F}_{\lambda}$ that the first $V(\lambda)$ $p$-values correspond to the true null and $p_{1},\ldots,p_{V(\lambda)}\leq \lambda$. In particular, we may consider $p^{(1)}=(0,p_{2},p_{3},\ldots,p_{m})$ and $p^{(2)}=(0,0,p_{3},\ldots,p_{m})$  if $V(\lambda)\geq 1$ and $V(\lambda)\geq 2$, respectively. Note that we introduced $p^{(j)}$ for $j>1$ in \Cref{sec:higher_moments}. Since $\widehat \alpha_{R:m}\leq \lambda$ we deduce from (BI\ref{BI3}) that 
\begin{align*}
		E\Bigl( \Bigl(\frac{V}{R}\Bigr)^2\Bigl | {\mathcal F }_{\lambda} \Bigr)&=V(\lambda) \E \Bigl( \frac{ \mathbf{1}\{p_1\leq \widehat \alpha_{R:m}\} }{ R^2} \Bigl| {\mathcal F }_{\lambda}\Bigr) \\ 
		&+ V(\lambda) (V(\lambda)-1)\E \Bigl( \frac{\mathbf{1}\{p_{1}\leq \widehat \alpha_{R:m},p_{2}\leq \widehat \alpha_{R:m}\}}{ R^2} \Bigl |{\mathcal F }_{\lambda}\Bigr).
	\end{align*} 
	Note that $p_{1},\ldots,p_{V(\lambda)}$ conditioned on ${\mathcal F }_{\lambda}$ are i.i.d. uniformly distributed on $(0,\lambda)$ if $V(\lambda)>0$. It is easy to see that $p_{1}\leq \widehat \alpha_{R:m}$ implies $R=R^{(1,\lambda)}$ and we have $P(p_{1}\in(\widehat\alpha_{R:m},\widehat\alpha_{R^{(1,\lambda)}:m}])=0$. Both were already known and used, for instance, in Heesen and Janssen \cite{HeesenJanssen2015,HeesenJanssen2016}. Since $p_{1}$ and $R^{(1,\lambda)}$ are independent conditionally on ${\mathcal F }_{\lambda}$ we obtain from Fubini's Theorem that
	\begin{align*} 
		\E \Bigl( \frac{ \mathbf{1}\{p_{1}\leq \widehat \alpha_{R:m}\} }{ R^2} \Bigl| {\mathcal F }_{\lambda}\Bigr) = \E \Bigl( \frac{ \mathbf{1}\{p_{1}\leq \widehat \alpha_{R^{(1,\lambda)}:m}\} }{ (R^{(1,\lambda)})^2} \Bigl| {\mathcal F }_{\lambda}\Bigr)=\frac{\alpha}{\lambda\widehat m_0}\E \Bigl( \frac{ 1 }{ R^{(1,\lambda)}}  \Bigl| {\mathcal F }_{\lambda}\Bigr).
	\end{align*}         
	Hence, we get the second summand of the right-hand side in \eqref{enu:theo:form_moments_second_moment}. 
	To obtain the first term, it is sufficient to consider $V(\lambda)\geq 2$. Since $p_{1}$, $p_{2}$ and $R^{(2,\lambda)}$ are independent conditionally on ${\mathcal F }_{\lambda}$ we get  similarly to the previous calculation:
	\begin{align}\label{eqn:proof_second_moment_R2}
		&E \Bigl( \frac{\mathbf{1}\{p_{1}\leq \widehat \alpha_{R:m},p_{2}\leq \widehat \alpha_{R:m}\}}{ R^2} \Bigl |{\mathcal F }_{\lambda}\Bigr)\\
		&= E \Bigl( \frac{\mathbf{1}\{p_{1}\leq \widehat \alpha_{R^{(2,\lambda)}:m},p_{2}\leq \widehat \alpha_{R^{(2,\lambda)}:m}\}}{ (R^{(2,\lambda)})^2} \Bigl |{\mathcal F }_{\lambda}\Bigr)= \frac{\alpha^2}{\lambda^2}\frac{1}{\widehat m_0^2},
	\end{align}   
	 which completes the proof of \eqref{enu:theo:form_moments_second_moment}. Combining \eqref{enu:theo:form_moments_second_moment}, \eqref{eqn:FDR_forumla} and the variance formula $\text{Var}(Z)=\E(Z^2)-\E(Z)^2$ yields \eqref{enu:theo:form_moments_variance}. The proof of \eqref{enu:theo:form_moments_E(V)} is based on the same techniques as the one of \eqref{enu:theo:form_moments_E(V)}, to be more specific we have
 	\begin{align}\label{eqn:E(V|F)}
		\E(V|\mathcal{F}_{\lambda,m}) = V(\lambda) P( p_{1}\leq \widehat \alpha_{R^{(1,\lambda)}:m}  \:|\: {\mathcal F }_{\lambda} ) = V(\lambda)\E\Bigl( \frac{R^{(1,\lambda)}\alpha}{\lambda\widehat m_0}\Bigl | {\mathcal F }_{\lambda} \Bigr).
	\end{align}
		
\subsection{Proof of \Cref{lem:bound_var}} To improve the readability of the proof, all indices $m$ are submitted except for $K_m$.

    \underline{\eqref{enu:lem_bound_var_Cm}:}  By \Cref{theo:form_moments}\eqref{enu:theo:form_moments_variance} and (A\ref{enu:ass_A_mo_leq_R}) it remains to show that
	\begin{align*}             
		\Bigl( \frac{\alpha}{\lambda} \Bigr)^2\E \Bigl( \frac{V(\lambda)}{\widehat m_0^2} \Bigr) \leq \E \Bigl( \frac{V(\lambda)}{R(\lambda)^2} \Bigr) \leq \E \Bigl( \frac{\mathbf{1}\{V(\lambda)>0\}}{V(\lambda)} \Bigr)
	\end{align*}
	is smaller than $2/(\lambda(m_0+1))$. It is known and can easily be verified that $\E(\mathbf{1}\{X>0\} X^{-1})\leq 2\E((1+X)^{-1})\leq 2 p^{-1}(n+1)^{-1}$ for any Binomial-distributed $X\sim B(n,p)$. Since $V(\lambda)\sim B(m_0,\lambda)$ we obtain the desired upper bound, see also p. 47ff of \citet{HeesenJanssen2016} for details.
		
	\underline{\eqref{enu:lem_bound_var_E(ind_E(V|F))}:} We can deduce from \eqref{eqn:E(V|F)} that
	\begin{align*}                                
		Y_{\lambda}:= \frac{\mathbf{1}\{V(\lambda)>0\}}{\E(V|\mathcal{F}_{\lambda})} = \frac{\lambda}{\alpha}\frac{\widehat m_0}{V(\lambda)} \frac{\mathbf{1}\{V(\lambda)\geq 1\} }{\E( R^{(1,\lambda)}|\mathcal{F}_{\lambda})}.
	\end{align*}   
	Note that
	\begin{align*}
		P\Bigl( \E(V|\mathcal{F}_{\lambda} \leq t\Bigr)\leq P(V=0)+P\Bigl( Y_{\lambda}\geq \frac{1}{t} \Bigr).
	\end{align*}                                                                                                              
	Thus, by Markoff's inequality it remains to verify $\E(Y_{\lambda})\leq D_{\lambda}$. We divide the discussion of $\E(Y_{\lambda})$ into two parts. We obtain from $R^{(1,\lambda)}\geq 1$ and Hoeffding's inequality, see p. 440 in \citet{ShorackWellner2009}, that
	\begin{align*}
		\E \Bigl( Y_{\lambda}\mathbf{1}\Bigl\{ \frac{V(\lambda)}{m_0}\leq \frac{\lambda}{2}\Bigr\} \Bigr) \leq \frac{\lambda m_0K_m}{\alpha} P\Bigl( V(\lambda)\leq \frac{\lambda m_0}{2} \Bigr) \leq \frac{\lambda m_0 K_m}{\alpha}\exp\Bigl(-\frac{1}{8}m_0\lambda\Bigr).
	\end{align*}
	Moreover, we obtain from  Jensen's inequality and \Cref{theo:form_moments}\eqref{enu:theo:form_moments_variance} that    
	\begin{align*}
		&\E \Bigl( Y_{\lambda}\mathbf{1}\Bigl\{ \frac{V(\lambda)}{m_0}\geq \frac{\lambda}{2}\Bigr\} \Bigr) \leq \frac{\lambda}{\alpha}\Bigl(\frac{2K_m}{\lambda} \Bigr)^2\E \Bigl( \frac{V(\lambda)}{\widehat m_0} \E\Bigl(\frac{1}{R^{(1,\lambda)}}|\mathcal{F}_{\lambda}\Bigr) \Bigr)\\
		&\leq  \Bigl( \frac{2K_m}{\lambda} \Bigr)^2\Bigl[\frac{\lambda^2}{\alpha^2}\text{Var}\Bigl( \frac{V}{R} \Bigr) +  \E \Bigl( \frac{V(\lambda)}{\widehat m_0^2} \Bigr) -\text{Var}\Bigl( \frac{V(\lambda)}{\widehat m_0} \Bigr)\Bigr].
	\end{align*}
	Finally, combining this with \eqref{eqn:E(V/m0^2)} yields the statement.
		
	\subsection{Proof of \Cref{theo:consistency}}
	Since $V_m/R_m$ is bounded by $1$ the consistency statement in \eqref{enu:theo:consis_consis} is equivalent to $\text{Var}(V_m/R_m)\to 0$. Due to $V_m(\lambda)/m\to \kappa_0\lambda$ a.s. and $K\geq (\widehat m_0/m)\geq (\alpha/\lambda)(V_m(\lambda)/m)\to \alpha\kappa_0$ a.s. we deduce from \Cref{lem:bound_var} that  $\text{Var}(\text{FDP}_m)\to 0$ is equivalent to $\text{Var}(V_m(\lambda)/\widehat m_0)\to 0$ and $\E((R_m^{(1,\lambda)})^{-1}|\mathcal{F}_{\lambda,m})\to 0$ in $P_m$-probability. Moreover, we can conclude from this that $\text{Var}(V_m(\lambda)/\widehat m_0)\to 0$ is equivalent to \eqref{eqn:theo:consistency_m0/m-E()}. Since $R_m^{(1,\lambda)}\geq 1$ we have equivalence of $\E((R_m^{(1,\lambda)})^{-1}|\mathcal{F}_{\lambda,m})\to 0$ in $P_m$-probability and \eqref{eqn:theo:consistency_R1}.
	
	\subsection{Proof of \Cref{theo:Vm_to_infty}}
	As a consequence of the assumptions we have at least for a subsequence $n(m)\to \infty$ that  
	\begin{align*}
		\lim_{m\to\infty}\E\Bigl(\frac{V_{n(m)}(\lambda)}{\widehat m_0}\Bigr)\to C\in \Bigl[ \frac{1}{K},\frac{\lambda}{\alpha} \Bigr].
	\end{align*}
	We suppose, contrary to our claim, that $V_{m,\alpha}$ does not converge to $\infty$ in $P_m$-probability for some $\alpha\in(\alpha_1,\alpha_2)$. Since $\alpha\mapsto V_{m,\alpha}$ is increasing we can suppose without loss of generality that $\lambda^{-1}\alpha C\notin \mathbb{Q}$ (otherwise take a smaller $\alpha>\alpha_1$). By our contradiction assumption there is some $k\in\N\cup\{0\}$ and a subsequence of $\{n(m):m\in\N\}$, which we denote by simplicity also by $n(m)$, with $n(m)\to\infty$ such that $P_{n(m)}(V_{n(m),\alpha}=k)\to \beta\in(0,1] $. We can deduce from \eqref{eqn:FDR_forumla} and the consistency that 
	\begin{align*}
		(V_{n(m),\alpha}/R_{n(m),\alpha}) \mathbf{1}\{V_{n(m),\alpha}=k\} -(\alpha/\lambda) C  \mathbf{1}\{V_{n(m),\alpha}=k\} \to 0
	\end{align*}
	in  $P_{n(m)}$-probability. In particular, it holds that
	\begin{align*}
		P_{n(m)}\Bigl(R_{n(m),\alpha}=\frac{k\lambda}{C\alpha},V_{n(m),\alpha}=k\Bigr) \to\beta>0,
	\end{align*}
	which leads to a contradiction since $(\lambda k)/(C\alpha)\notin \N\cup\{0\}$.
		
\subsection{Proof of \Cref{exam:consistency}}	
	Clearly, the $p$-values from the false null are i.i.d. with distribution function $F_1$ given by $F_1(t)=t\mathbf{1}\{t<x_0\}+\mathbf{1}\{t\geq x_0\}$. From \Cref{theo:cons_kappa=1} with $t_m=x_0$ regarding \Cref{rem:cons_kappa=1}\eqref{enu:rem:cons_kappa=1_F1} and straight forward calculations we obtain the consistency of BH$(2\alpha)$ and Stor$(\alpha,\lambda)$, for the latter see also \Cref{cor:Storey}\eqref{enu:cor:Storey_consis}. \\
	Let us now have a look at BH$(\alpha)$. First, we compare the empirical distribution function $\widehat F_m$ and the Simes line $t\mapsto g(t):=t/\alpha=4t$. From the Glivenko-Cantelli Theorem it is easy to see that $\widehat F_m$ tends uniformly to $F$ given by $F(t)= t\mathbf{1}\{t< x_0\}+ 1/2(t+1)\mathbf{1}\{t\geq x_0\} $. Clearly, the Simes line lies strictly above $F$ on $(0,1)$, whereas the Simes line $t\mapsto t/(2\alpha)$ corresponding to BH$(2\alpha)$ hits $F$.  Hence, $P_m(\sup_{\varepsilon\leq x \leq 1}\widehat F_m(x)-g(x)<0)\to 1$ for all $\varepsilon>0$. Let $0<\widetilde \lambda<x_0$. Then $P_m(\alpha_{R_m:m}\leq \widetilde \lambda)\to 1$ follows. That is why we can restrict our asymptotic considerations to the portion of $p$-values with $p_{i,m}\leq \widetilde \lambda$ and the inconsistency follows, compare to the proof of \Cref{theo:nec_for_cons}\eqref{enu:theo:nec_for_cons_unif_(0,lamb)}.
	
\subsection{Proof of \Cref{theo:nec_for_cons}}
	\underline{\eqref{enu:theo:nec_for_cons_m1_to_infty}:} We suppose contrary to our claim that $\liminf_{m\to\infty}m_1=k\in \N\cup\{0\}$. Then $m_1=k$ for infinitely many $m\in\N$. Turning to subsequences we can assume without a loss that $m_1=k$ for all $m\in\N$. Note that  (A\ref{enu:ass_A_kappa0}) holds for $\kappa_0=1$ in this case. Hence, it is easy to see that
	$\liminf_{m\to\infty}R_m(\lambda)/m\geq \lambda \text{ a.s.}$
	Combining this and (A\ref{enu:ass_A_mo_geq_CR}) yields
	\begin{align*}\label{eqn:proof_nec_hatalpha<tildealpaha}
		P_m\Bigl( \widehat \alpha_{i:m}\leq \Bigl( \frac{i}{m}\widetilde\alpha \Bigr)\text{ for all }i=1,\ldots,m \Bigr)\to 1 \text{ with }\widetilde \alpha=\frac{1}{C}<1.
	\end{align*}
	Hence, we can deduce from \Cref{exam:intro_extreme_cases}\eqref{enu:exam:intro_extreme_cases_DU} that
	\begin{align*}
		\liminf_{m\to\infty} P(V_m=0) \geq \liminf_{m\to\infty}P_m(V_m^{\text{BH}}(\widetilde \alpha,k)=0)>0.
	\end{align*}
	But this contradicts the necessary condition \eqref{eqn:nec_cond_cons_P(V>0)} for consistency.\\
		
	\underline{\eqref{enu:theo:nec_for_cons_unif_(0,lamb)}:} Suppose for a moment that we condition on $\mathcal{F}_{\lambda,m}$. Hence, $R_m(\lambda)$ and $\widehat m_0$ can be treated as fixed numbers. Without loss of generality we assume that $p_{1,m},\ldots,p_{R_m(\lambda),m}\leq \lambda$. Define new $p$-values $q_{1,R_m(\lambda)},\ldots,q_{R_m(\lambda),R_m(\lambda)}$ by $q_{i,R_m(\lambda)}:=p_{i,m}/\lambda$ for all $i=1,\ldots,R_m(\lambda)$. The values $V_m$ and $R_m$ are the same for the step-up test for $(p_{i,m})_{i\leq m}$ with critical values $\widehat \alpha_{i:m}=(i/\widehat m_0)\alpha$ and for $(q_{i,R_m(\lambda)})_{i\leq R_m(\lambda)}$ with critical values $\widehat \alpha^{(q)}_{i:R_m(\lambda)}=(i/R_m(\lambda))\widetilde \alpha_m$, where $\widetilde \alpha_m=(R_m(\lambda)/\widehat m_0)(\alpha/\lambda)$. In the situation here $q_{1,R_m(\lambda)},\ldots,q_{R_m(\lambda),R_m(\lambda)}$ are i.i.d. uniformly distributed on $(0,1)$ and so correspond to a DU$(R_m(\lambda),0)$-configuration. 
	That is why 
	\begin{align}\label{eqn:P(Vm=0|F)}
		P_m(V_m=0|\mathcal{F}_{\lambda,m}) = P_m(V_{R_m(\lambda)}^{\text{BH}}(\widetilde \alpha_m,0)=0|\mathcal{F}_{\lambda,m}).
	\end{align}
	Since by the strong law of large numbers and  (A\ref{enu:ass_A_mo_geq_CR}) we have $R_m(\lambda)\to \infty$ a.s. and $P_m(\widetilde \alpha_m \leq C^{-1})\to 1$ we can conclude from  \eqref{eqn:P(Vm=0|F)} and \Cref{exam:intro_extreme_cases}\eqref{enu:exam:intro_extreme_cases_DU}
	\begin{align*}
		\liminf_{m\to\infty} P_m(V_m=0) \geq (1-C^{-1})
	\end{align*}
	and so the necessary condition \eqref{eqn:nec_cond_cons_P(V>0)} for consistency is not fulfilled.
		
	%Observe that $\frac{R_m(\lambda)}{m}\geq \frac{m_1}{m}\lambda+\frac{m_0}{m}\frac{V_m(\lambda)}{m_0}\to \lambda \text{ a.s.}$ and so \eqref{eqn:proof_nec_hatalpha<tildealpaha} holds again. 
	
	\subsection{Proof of \Cref{lem:R_m/m_to_0}}
	Analogously to the proof of \Cref{theo:nec_for_cons}\eqref{enu:theo:nec_for_cons_unif_(0,lamb)} we condition under $\mathcal{F}_{\lambda,m}$ and introduce the new $p$-value $q_{i,R_m(\lambda)}$ and the new critical value $\widehat \alpha^{(q)}_{i:R_m(\lambda)}$ for $i\leq R_m(\lambda)$ as well as the new level $\widetilde\alpha_m$. The respective empirical distribution functions of the new $p$-value $(q_{i,m})_{i\leq R_m(\lambda)}$ are denoted by $\widehat F_{R_m(\lambda)}^{(q)}, \widehat F_{0,R_m(\lambda)}^{(q)}, \widehat F_{1,R_m(\lambda)}^{(q)}$, compare to the definition of $\widehat F_{j,m}$ in \Cref{theo:cons_kappa=1}. Note that $\widehat \alpha_{R_m:R_m(\lambda)}^{(q)}$ is the largest intersection point of $\widehat F_m^{(q)}$ and the Simes line $t\mapsto f_{\widetilde \alpha_m}(t)= :t/\widetilde\alpha_m$. Note that  $R_m(\lambda)\geq m_0 \widehat F_{0,m}(\lambda)\to \infty$ $P_m$-a.s. and, hence, by (A\ref{enu:ass_A_mo_geq_CR}) $P_m(\widetilde \alpha_m\leq C^{-1})\to 1$. From this and the Glivenko-Cantelli Theorem  we obtain that for all $\varepsilon\in(0,1)$
	\begin{align*}
		P_m\Bigl( \,\sup_{t\in[\varepsilon,1]} \widehat F_{0,R_m(\lambda)}^{(q)}(t) - f_{\widetilde \alpha_m}(t)\leq (1-C)\varepsilon\Bigl | \mathcal{F}_{\lambda,m} \Bigr) \to 1.
	\end{align*}
	Combining this and $\widehat F_{R_m(\lambda)}^{(q)}(t)-\widehat F_{0,R_m(\lambda)}^{(q)}(t)\leq m_1/m\to 0$ we can deduce that $P_m(\widehat\alpha_{R_m:m}\leq \lambda \varepsilon)\to 1$ for all $\varepsilon>0$. 
		
	\subsection{Proof of \Cref{theo:cons_kappa=1}}
		Clearly, all $p_{i,m}\leq t_m$ are rejected and, in particular, $V_m\geq V_m(t_m)$ if $p_{R_m(t_m):m}\leq (R_m(t_m)/\widehat m_0)\alpha$. The latter is fulfilled if $(t_m/\alpha) \widehat m_0 \leq R_m(t_m)$, or equivalently
	\begin{align}\label{eqn:cons_ka=1_suff}
		\frac{\widehat m_0}{m\alpha} \leq \frac{m_0}{m} \frac{\widehat F_{0,m}(t_m)}{t_m}+\frac{m_1}{m} \frac{\widehat F_{1,m}(t_m)}{t_m}.
	\end{align}
	Note that by Chebyshev's inequality
	\begin{align*}
		P_m\Bigl( \frac{m_0}{m} \frac{\widehat F_{0,m}(t_m)}{t_m} \geq \kappa_0 - \frac{1}{2}\delta \Bigr) \geq 1 - \frac{1}{mt_m}\Bigl(\frac{1}{2}\delta+\frac{m_0}{m}-\kappa_0\Bigr)^{-2}\to 1.
	\end{align*}
	Combining this, \eqref{eqn:theo:consistency_m0/m}, \eqref{eqn:cons_kappa=1_condition} and \eqref{eqn:cons_ka=1_suff} yields
	\begin{align*}
		P_m(V_m\geq V_m( t_m)) \to 1.
	\end{align*}
	Since $V_m(t_m)\sim B(m_0,t_m)$ and $m_0t_m\to \infty$ the statement follows. 
		
	\subsection{Proof of \Cref{rem:cons_kappa=1}}
	By \Cref{theo:cons_kappa=1} it remains to show that
	\begin{align*}
		&P_m\Bigl( \frac{m_1}{m}\frac{ \widehat F_{1,m_1}(t_m) }{t_m} \geq \frac{1}{2}\delta - \kappa_0 + \frac{1}{\alpha}C_0  \Bigr) \nonumber \\
		&=P_m\Bigl( \sqrt{m_1} \frac{ \widehat F_{1,m_1}(t_m)-F_1(t_m)}{\sqrt{F_1(t_m)(1-F_1(t_m))}} \geq \sqrt{m_1}t_m \frac{ \frac{m}{m_1}(\frac{\delta}{2}-\kappa_0+\frac{1}{\alpha}C_0 - \frac{m_1}{m}\frac{F_1(t_m)}{t_m})}{\sqrt{F_1(t_m)(1-F_1(t_m))}}\Bigr) %\label{eqn:rem:cons_ka=1_suff}
	\end{align*}
	converges to $1$. Note that the left-hand side of the last row converges in distribution to $Z\sim N(0,1)$. Moreover, by straightforward calculations it can be concluded from \eqref{eqn:cons_kappa=1_cond_F1} and $C_0\geq \kappa_0\alpha$ that the right-hand side tends to $-\infty$, which completes the proof.
		
\subsection{Proof of \Cref{lem:comp_BH_adapt_consis}}
	It is easy to see that \eqref{eqn:lem_comp_BH_condition_to_1} always holds if $C_0<1$. From \eqref{eqn:lem_comp_BH_condition_to_1} we obtain immediately that
	\begin{align*}
		&P_m\Bigl( \max_{i=1,\ldots,m}\Bigl\{ \widehat \alpha^{\text{BH}}_{i,m}-\widehat \alpha_{i,m} \Bigr\} \leq 0 \Bigr)\leq P_m\Bigl( \frac{1}{m}-\frac{1}{\widehat m_0} \leq 0 \Bigr) \to 1\\
		& \text{ and so }P_m\Bigl( R_m^{(1,\lambda)}\geq R_m^{(1,\lambda),\text{BH}}  \Bigr)\to 1,
	\end{align*}
	where $R_m^{(1,\lambda),\text{BH}}$ is the corresponding random variable for the BH test. Now, suppose that we have consistency for the BH test. Then combining \Cref{theo:consistency} with the above yields that $R_m^{(1,\lambda),\text{BH}}$ and so $R_m^{(1,\lambda)}$ converges to infinity in $P_m$-probability. Finally, we deduce the consistency of the adaptive test from \Cref{theo:consistency}.
		
\subsection{Proof of \Cref{cor:Storey}}
	\underline{\eqref{enu:cor:Storey_ka=1}:} Clearly, $\widetilde m_0(\lambda_{i-1},\lambda_i)/m\to 1$ a.s. for all $i=1,\ldots,k$ and $R_m(\lambda)/m\to \lambda$ a.s. Thus, \eqref{eqn:theo:consistency_m0/m} holds for $C_0=1$. Finally, \eqref{eqn:cor:Storey_ka=1_FDR} follows from \eqref{eqn:FDR_forumla}.
		
	\underline{\eqref{enu:cor:Storey_ka<1}:}
	First, we introduce new estimators $\widetilde m_{0,i}$ and new weights $\widehat \gamma_{i,m}\geq 0$ for all $i=1,\ldots,k$:
	\begin{align*}
		\widetilde m_{0,i}:=m\frac{ 1-\widehat F_m(\lambda_{i-1})-\frac{i}{m} }{1-\lambda_{i-1}} \text{ and }\widehat \gamma_i:= \Bigl( \frac{\widehat\beta_{i,m}}{\lambda_i-\lambda_{i-1}} - \frac{\widehat\beta_{i-1,m}}{\lambda_{i-1}-\lambda_{i-2}} \Bigr)(1-\lambda_{i-1}),
	\end{align*}
	where we use $\widehat \beta_{0,m}:=0$. It is easy to check $\widetilde m_0=\sum_{i=1}^k \widehat \gamma_{i,m} \widetilde m_{0,i}$ and $\sum_{i=1}^k \widehat \gamma_{i,m} =1$. From \eqref{eqn:cor:Storey_condition_liminf} and the strong law of large numbers it follows 
	\begin{align}\label{eqn:proof_Storey_as_V_m0}
		\frac{V_m(\lambda)}{m} \to\kappa_0\lambda\text{ a.s. and }\limsup_{m\to\infty}\frac{\widetilde m_{0,i}}{m}\leq 1-\frac{\varepsilon_i}{1-\lambda_i} \text{ a.s.}
	\end{align}
	In particular, by \eqref{eqn:cor:Storey_condition_eps>0}
	\begin{align*}
		\limsup_{m\to\infty}\frac{\widetilde m_0}{m}\leq  1- \frac{\delta(1-\lambda_{j-1})\varepsilon_j}{1-\lambda_j}\leq \frac{1}{1+\delta_0}\;\text{ a.s.}
	\end{align*}
	for some $\delta_0>0$. Consequently, 
	\begin{align*}
		\liminf_{m\to\infty} \frac{V_m(\lambda)}{\widetilde m_0}\geq \lambda \kappa_0(1+\delta_0)\text{ a.s.}
	\end{align*}
	It is easy to verify that (A\ref{enu:ass_comb_Storey_mo_geq_al_lam_R}) implies $P_m(\widetilde m_0(\lambda_{i-1},\lambda_{i})> C_i(\alpha/\lambda)R_m(\lambda))\to 1$ for appropriate $C_i>1$ and for all $i$. Hence, (A\ref{enu:ass_A_mo_geq_CR}) is fulfilled and, in particular, $P_m(\widehat m_0 = \widetilde m_0)\to 1$. Finally, we obtain the statement from \eqref{eqn:FDR_forumla}.
		
	\underline{\eqref{enu:cor:Storey_consis}:} Define $\widetilde m_{0,i}$ and $\widehat \gamma_{i,m}$ as in the proof of \eqref{enu:cor:Storey_ka<1}. Then,
	\begin{align*}
		\widehat{\gamma}_i\to\frac{\beta_{i}}{\lambda_{i}-\lambda_{i-1}}- \frac{\beta_{i-1} }{\lambda_{i-1}-\lambda_{i-2}}=:\gamma_i\text{ a.s.}
	\end{align*}
	for all $i=1,\ldots,k$. Clearly, (A\ref{enu:ass_A_mo_leq_Km}) and (A\ref{enu:ass_A_mo_geq_CR}) are fulfilled, see for the latter the end of the proof of \eqref{enu:cor:Storey_ka<1}. Moreover, \eqref{eqn:theo:consistency_m0/m} holds for some $C_0\in[0,1]$ since
	\begin{align*}
		\frac{\widehat m_0}{m}\to 1 - (1-\kappa_0)\sum_{i=1}^k \frac{\gamma_i\varepsilon_{i}}{1-\lambda_i}\text{ a.s.}
	\end{align*}
	Due to \eqref{eqn:cor:Storey_condition_eps>0_consis} we have $C_0<1$ iff $\kappa_0<1$. Consequently, by \Cref{theo:cons_kappa=1} and \Cref{lem:comp_BH_adapt_consis} it remains to verify \eqref{eqn:lem_comp_BH_condition_to_1}  in the case of $\kappa_0=1$. \\
	Consider $\kappa_0=1$. By assumption we have $m_1/\sqrt{m}\to\infty$ in this case. First, observe that by the central limit theorem it holds for all $i=1,\ldots,k$ that
	\begin{align}\label{eqn:cor:Storey_Zni}
		Z_{i,m}:=\sqrt{m} \frac{m_0}{m} \Bigl( \frac{1-\widehat F_{0,m}(\lambda_i)}{1-\lambda_i}-1 \Bigr) \overset{\mathrm d}{\longrightarrow} Z_i\sim N(0,\sigma_i^2)
	\end{align}
	for some $\sigma_i\in(0,\infty)$. Let $\xi:=\varepsilon_j/(8(1-\lambda_j))>0$. By \eqref{eqn:cor:Storey_condition_lim} and  \eqref{eqn:cor:Storey_condition_eps>0_consis} 
	\begin{align}
		&P_m\Bigl( \frac{1-\widehat F_{1,m}(\lambda_j)}{1-\lambda_j}\leq 1 - 4\xi \Bigr)\to 1 \label{eqn:cor_proof_Storey_1-F1_j}\\
		\text{ and }&P_m\Bigl( \frac{1-\widehat F_{1,m}(\lambda_i)}{1-\lambda_i}\leq 1 + \frac{\xi}{2}\gamma_j \Bigr)\to 1. \label{eqn:cor_proof_Storey_1-F1_i}
	\end{align}
	for all $i\in\{1,\ldots,k\}\setminus\{j\}$. Moreover, from \eqref{eqn:cor:Storey_Zni} we get
	\begin{align*}
		&P_m\Bigl( \frac{m_0}{m}\frac{1-\widehat F_{1,m}(\lambda_j)+\frac{1}{m}}{1-\lambda_j}+\frac{m_1}{m}\Bigl(  1 - 4\xi + \frac{\frac{1}{m}}{1-\lambda_j} \Bigr)\leq 1-\frac{m_1}{m} 2\xi \Bigr)\\
		&= P_m\Bigl( Z_{i,m} \leq \frac{m_1}{\sqrt{m}}2\xi -\frac{1}{\sqrt{m}}\frac{1}{1-\lambda_{i}}\Bigr)\to 1.
	\end{align*}
	From this and \eqref{eqn:cor_proof_Storey_1-F1_j} $P_m(\widetilde m_{0,j}\leq  1-(m_1/m) 2\xi)\to 1$ follows. Analogously, we obtain from \eqref{eqn:cor_proof_Storey_1-F1_i} that $P_m(\widetilde m_{0,i}\leq  1+(m_1/m) \gamma_j\xi)\to 1$ for all $i\neq j$.  Since $\sum_{i=1, i\neq j}^{k}\widehat \gamma_i\leq 1$ and $P_m(2\widehat \gamma_j \geq \gamma_j)\to 1$ we can finally conclude \eqref{eqn:lem_comp_BH_condition_to_1}.
		
	% that $P_m(R_m<r_m)\to 0$ with $r_m=\lfloor \log(m)\rfloor=\max\{k\in\N:k\leq \log(m)\}$. For this purpose define a new triangular schema $(q_{i,m_0})_{i\leq m_0}$ of $p$-value consisting only of the true $p$-value, i.e. $q_{i,m_0}=p_{i,m}$ for all $i=1,\ldots,m_0$. Then
		%\begin{align*}
		%	P_m(R_m<r_m) &\leq P_m( q_{r_m:m_0}\geq \widehat \alpha_{r_m:m_0}) \\
		%	&= P_m\Bigl( Z_m\geq r_m \Bigr)\text{ with }Z_m:=\sum_{i=1}^{m_0}\mathbf{1}\Bigl\{ q_{i,m_0}\leq\widehat \alpha_{r_m:m_0} \Bigr\}.
		%\end{align*} 
		%Note that conditioned on $\mathcal{F}_{\lambda,m}$ the distribution of $Z_m$ is $B(V_m(\lambda),\widehat \alpha_{r_m:m_0})$.
\subsection{Proof of \Cref{lem:cons_gen_alp}} 
\underline{\eqref{enu:lem:cons_gen_alp_suff_cond}:} First, we introduce for $j=1,2$:
\begin{align}\label{eqn:lem_cons_gen_alp_defi_proof_psi}
	\psi_{m,j}:=\frac{m_0\alpha}{m+b_m-a_mR_m^{(j,1)}}.
\end{align}
Using the formulas presented at the beginning of \Cref{sec:AORC} we obtain:
\begin{align*}
	\text{Var}\Bigl( \frac{V_m}{R_m} \Bigr)= &\E \Bigl( \frac{1}{R_m^{(1,1)}} \psi_{m,1} \Bigr)  + \text{Var}(\psi_{m,1}) + \E\Bigl( \psi_{m,2}^2-\psi_{m,1}^2 \Bigr)- \frac{1}{m_0}\E( \psi_{m,2} )^2
	% &+\E\Bigl( \frac{a_m( R_m^{(2,1)}-R_m^{(1,1)})}{\alpha m_0} [\psi_{m,1}\psi_{m,2}^2+\psi_{m,1}^2\psi_{m,2}] \Bigr).
\end{align*}
Note that by $R_m^{(j,1)}\leq m$, $0\leq b_m\leq m$ for large $m$ and \eqref{eqn:alpha_im_quotient} we get: 
\begin{align}\label{eqn:lem_cons_gen_alp_psi_tight}
	\kappa_0\leftarrow\frac{m_0}{m+b_m}\leq \psi_{m,j}\leq \frac{m_0}{\alpha m}\to \frac{\kappa_0}{\alpha}.
\end{align}
Hence, the fourth summand $-\E( \psi_{m,2} )^2/m_0$ in the formula for \text{Var}$(V_m/R_m)$  tends always to $0$. Since, clearly, the first three summands are non-negative it remains to show that each of these summands tends to $0$ iff our conditions  \eqref{eqn:lem_cons_gen_alp_R1_to_infty}-\eqref{eqn:lem:cons_gen_alp_def_psi1} are fulfilled. By \eqref{eqn:lem_cons_gen_alp_psi_tight} we have equivalence of \eqref{eqn:lem:cons_gen_alp_def_psi1} and Var$(\psi_{m,1})\to 0$, as well as of \eqref{eqn:lem_cons_gen_alp_R1_to_infty} and $\E ( \psi_{m,1}/R_m^{(1,1)}  )$ to $0$. Observe that
$\psi_{m,2}-\psi_{m,1} =  Z_m\psi_{m,1}\psi_{m,2}$ with $Z_m:=(a_m/m_0)(R_m^{(2,1)}-R_m^{(1,1)})$. From \eqref{eqn:lem_cons_gen_alp_psi_tight} and $0\leq Z_m\leq  1-\alpha$ we obtain  that  $E(\psi_{m,2}-\psi_{m,1})\to 0$ iff \eqref{eqn:lem_cons_gen_alp_R2-R1} holds. Finally, combining this, \eqref{eqn:lem_cons_gen_alp_psi_tight}, $\psi_{m,2}^2-\psi_{m,1}^2=(\psi_{m,2}-\psi_{m,1})(\psi_{m,2}+\psi_{m,1})$ and $\psi_{m,2}\geq \psi_{m,1}$ yields that $\E( \psi_{m,2}^2-\psi_{m,1}^2 )\to 0$ iff \eqref{eqn:lem_cons_gen_alp_R2-R1} is fulfilled.

\underline{\eqref{enu:lem:cons_gen_alp_suff_cond_sparse}:} Similarly to \Cref{lem:R_m/m_to_0} we obtain by considering the (least favorable) DU$(m,m_1+j)$-configuration that $\alpha_{R_m^{(j,1)}:m}\to \infty$ and so $R_m^{(j,1)}/m\to \infty$ both in $P_m$-probability for $j=1,2$. Clearly, \eqref{eqn:lem_cons_gen_alp_R2-R1} follows. Moreover, we deduce from this, $a_m\leq 1-\alpha$ and $b_m/m\to 0$ that $\psi_{m,j}$ defined by \eqref{eqn:lem_cons_gen_alp_defi_proof_psi} converges to $\alpha$ in $P_m$-probability for $j=1,2$. This implies \eqref{eqn:lem:cons_gen_alp_def_psi1} and $E( V_m/R_m ) =  \E ( \psi_{m,1} )\to \alpha$. In particular, we have consistency by \eqref{enu:lem:cons_gen_alp_suff_cond}.

\subsection{Proof of \Cref{theo:LFC_ada_test}}
	\underline{\eqref{enu:theo:LFS_ada_test_cond_LFC_Exp}:} Let $p^*_i\in[0,1]$ be fixed for each $i\in I_{1,m}$. Let $P_m^*$ be the distribution fulfilling BI, where the $p_{i,m}\equiv p_i^*$ a.s. for all $i\in I_{1,m}$. From  \eqref{eqn:FDR_forumla} we get 
	\begin{align*}
		\int \frac{V_m}{R_m} \,\mathrm{ d }P_m^* = \frac{\alpha}{\lambda} \E\Bigl( \frac{V_m(\lambda)}{\widehat m_0} \Bigr).
	\end{align*}
	Moreover, we observe that the right-hand side only depends on $p_i^*$, $i\in I_{1,m}$, if $p_i^*>\lambda$. Consequently, we obtain the statement.
		
	\underline{\eqref{enu:theo:LFS_ada_test_cond_LFC_var}:} Due to \eqref{enu:theo:LFS_ada_test_cond_LFC_Exp} it remains to show that the conditional second moment is minimal under DU$_{\text{cond}}(m,M_{1,m}(\lambda))$. Clearly, BI and (A\ref{enu:ass_A_mo_leq_R}) are also fulfilled conditioned on $p^*_{\lambda,m}$. From \Cref{theo:form_moments}\eqref{enu:theo:form_moments_second_moment} we obtain
	\begin{align*}
		\E \Bigl( \Bigl( \frac{V_m}{R_m} \Bigr)^2 \Bigl | p^*_{\lambda,m} \Bigr)=  \E \Bigl( \frac{\alpha^2V_m(\lambda)(V_m(\lambda)-1)}{ \lambda^2\widehat m_0^2} + \frac{\alpha}{\lambda} \frac{ V_m(\lambda)}{\widehat m_0}  \frac{1}{R_m^{(1,\lambda)}}  \Bigl| p^*_{\lambda,m}\Bigr).
	\end{align*}
	It is easy to see that $V_m(\lambda)$ and $\widehat m_0$ are not affected and $R_m^{(1,\lambda)}$ increase if we set all $M_{1,m}(\lambda)$ $p$-value $p_{i,m}\leq \lambda$, $i\in I_{n,1}$, to $0$. 
		
	\underline{\eqref{enu:theo:LFS_ada_test_uncond_LFC_decrea}:} Since $V_m(\lambda)$ is not affected by any $p_{i,m}$, $i\in I_{1,m}$, the first statement follows from (A\ref{enu:ass_A_mo_increasing}) and \eqref{eqn:FDR_forumla}. If $p_{i,m}\leq \lambda$, $i\in I_{1,m}$,  decreases than $V_m(\lambda)$ and $\widehat m_0$ are not affected, and $R_m$ as well as $R_m^{(1,\lambda)}$ increase. Hence, the second statement follows from \Cref{theo:form_moments}\eqref{enu:theo:form_moments_variance}.
		
	\underline{\eqref{enu:theo:LFS_ada_test_uncond_LFC_max+min_DU}:} The statement follows immediately from \eqref{enu:theo:LFS_ada_test_uncond_LFC_decrea}.
	
\section{Appendix: Higher moments}\label{sec:higher_moments}
We extend the idea of the definition of $p_m^{(1,\lambda)}$ and $R_m^{(1,\lambda)}$ from \Cref{sec:consis}. For every $1<j\leq m_0$ we introduce a new $p$-value vector $p_{m}^{(j)}$ as a modification of $p_m=(p_{1,m},\ldots,p_{m,m})$ iteratively. If $V_m(\lambda)\geq j$ then we define $p_m^{(j,\lambda)}$ by setting $p_{i_k,m}$ equal to $0$ for $j$ different indices $i_1,\ldots,i_j\in I_{0,m}$ with $p_{i_k,m}\leq \lambda$, for convenience take the smallest $j$ indices with this property. Otherwise, if $V_m(\lambda)<j$ then set $p_m^{(j,\lambda)}$ equal to $p_m^{(j-1)}$. Moreover, let $R_m^{(j,\lambda)}=R_m^{(j,\lambda)}(p_m^{(j,\lambda)})$ be the number of rejections of the adaptive test for the (new) $p$-value vector $p_m^{(j,\lambda)}$. Note that $\widehat m_0$ is not affected by these replacements.
\begin{theorem}\label{theo:higher_moments}
	Under (A\ref{enu:ass_A_mo_leq_R}) we have for every $k\leq m$
	\begin{align*}
		&\E \Bigl( \Bigl( \frac{V_m}{R_m} \Bigr)^k  \Bigr)= \sum_{j=1}^{k} \alpha^j C_{j,k} \E\Bigl( \frac{ V_m(\lambda)\ldots(V_m(\lambda)-j+1)}{(\widehat m_0)^j} \E \Bigl( \Bigl( R^{(j,\lambda)}_m \Bigr)^{j-k}\Bigl | \mathcal{F}_{\lambda,m} \Bigr) \Bigr),\\
		&\text{where }C_{j,k}=\frac{1}{j!}\sum_{r=0}^{j-1}(-1)^{r}\binom{j}{r}(j-r)^k.
	\end{align*}
\end{theorem}
\begin{remark}\label{rem:higher_moments}
\begin{enumerate}[(a)]
		\item 	If we set $\widehat m_0=m_0$ and $\lambda=1$ then this formula coincide up to the factor $C_{j,k}$ with the result of \citet{FerreiraZwinderman2006}. By carefully reading their proof it can be seen that the coefficients $C_{r,k}$ have to be added. It is easy to check that $C_{1,k}=C_{k,k}=1$ but $C_{r,k}>1$ for all $1<r<k$. In particular, the coefficients $C_{j,2},\,C_{1,1}$, which are needed for the variance formula, are equal to $1$. 
				
		\item For treating one-sided null hypothesis the assumption (BI\ref{BI3}) need to be extended to i.i.d. $(p_{i,m})_{i\in I_{0,m}}$ $p$-values of the true null hypothesis, which are stochastically larger than the uniform distribution, i.e. $P(p_{i,m}\leq x)\leq x$ for all $x\in[0,1]$. In this case the equality in \Cref{theo:higher_moments} is not valid in general but the statement remains true if "$=$" is replaced by "$\leq$", analogously to the results of \citet{FerreiraZwinderman2006}.
\end{enumerate}
\end{remark}
\begin{proof}[Proof of \Cref{theo:higher_moments}]
	For the proof we extend the ideas of the proof of \Cref{theo:form_moments}. In particular, we condition on $\mathcal{F}_{\lambda,m}$. First, observe that
	\begin{align}\label{eqn:first_form_higher_mom}
		\E\Bigl( \frac{V_m^k}{R_m^k} \Bigl | \mathcal{F}_{\lambda,m}\Bigr) = \sum_{i_1,\ldots,i_k=1}^{V_m(\lambda)} \E \Bigl( \frac{\mathbf{1}\{p_{i_s,m}\leq \widehat{\alpha}_{R_m:m}, \, s\leq k\}}{R_m^{k}}\Bigl | \mathcal{F}_{\lambda,m}\Bigr).
	\end{align}
	Due to (BI\ref{BI3}) it is easy to see that each summand only depends on the number $j=\#\{i_1,\ldots,i_k\}$ of different indices. At the end of the proof we determine these summands in dependence of $j$.  But first we count the number of possibilities of choosing $(i_1,\ldots,i_k)$ which lead to the same $j$. Let $j\in\{1,\ldots,V_m(\lambda)\wedge k\}\}$ be fixed. Clearly, there are $\binom{V_m(\lambda)}{j}$ possibilities to draw $j$ different numbers $\{M_1,\ldots,M_j\}$ from the set $\{1,\ldots,V_m(\lambda)\}$. Moreover, by simple combinatorial considerations there are 
	\begin{align*}
		\sum_{r=0}^{j-1}(-1)^{r}\binom{j}{r}(j-r)^k
	\end{align*}
	possibilities of choosing indices $i_1,\ldots,i_k$ from $\{M_1,\ldots,M_j\}$ such that every $M_s$, $1\leq s \leq j$, is picked at least once, see e.g. (II.11.6) in \citet{Feller1968}. Consequently, we obtain from (BI\ref{BI3}) that \eqref{eqn:first_form_higher_mom} equals
	\begin{align*}
		\sum_{j=1}^{V_m(\lambda)\wedge k} C_{j,k}  V_m(\lambda)\ldots(V_m(\lambda)-j+1) \E \Bigl( \frac{\mathbf{1}\{p_{s,m}\leq \widehat{\alpha}_{R_m:m}, s\leq j\}}{R_m^k}\Bigl | \mathcal{F}_{\lambda,m} \Bigr). 
	\end{align*}
	Clearly, we can replace $V_m(\lambda)\wedge k$ by $k$ since each additional summand is equal to $0$. It remains to determine the summands. Let $j\leq V_m(\lambda)\wedge k$. Without loss of generality we can assume conditioned on $\mathcal{F}_{\lambda,m}$ that the first $V_m(\lambda)$ $p$-values correspond to the true null and $p_{1,m},\ldots,p_{V_m(\lambda),m}\leq \lambda$. In particular, we may consider $p_m^{(j,\lambda)}=(0,\ldots,0,p_{j+1,m},\ldots,p_{m,m})$. We obtain analogously to the calculation in \eqref{eqn:proof_second_moment_R2} and the one before it that
	\begin{align*}
		&\E \Bigl( \frac{\mathbf{1}\{p_{s,m}\leq \widehat{\alpha}_{R_m:m},s\leq j\}}{R_m^k}\Bigl | \mathcal{F}_{\lambda,m} \Bigr)\\
		&=	\E \Bigl( \frac{\mathbf{1}\{p_{s,m}\leq \widehat{\alpha}_{R^{(j,\lambda)}_m:m},s\leq j\}}{(R^{(j,\lambda)}_m)^k}\Bigl | \mathcal{F}_{\lambda,m} \Bigr)
		= \Bigl( \frac{\alpha}{\widehat m_0} \Bigr)^j \E\Bigl( \Bigl( R^{(j,\lambda)}_m \Bigr)^{j-k} \Bigl| \mathcal{F}_{\lambda,m} \Bigr).
	\end{align*} 
\end{proof}

\end{document}